 \documentclass[leqno]{article}
 \usepackage{amstex,theorem}
 \input epsf \epsfverbosetrue

 \newtheorem{Cor}{Corollary}[section] 
 \newtheorem{Prop}[Cor]{Proposition}
 \newtheorem{Thm}[Cor]{Theorem}
 \newtheorem{RThm}[Cor]{Residue Theorem}
 \newtheorem{DThm}[Cor]{Duality Theorem}
 \newtheorem{Con}[Cor]{Main Conclusion}
 \newtheorem{Prob}[Cor]{Problem}
 {
 \theorembodyfont{\rmfamily}
 \newtheorem{Def}[Cor]{Definition}
 \newtheorem{Ex}[Cor]{Example}
 \newtheorem{Exs}[Cor]{Examples}
 \newtheorem{Rem}[Cor]{Remark}
 }

 \newcommand{\no}{\refstepcounter{Cor}\medskip\noindent%
(\arabic{section}.\arabic{Cor})\enspace}

 \makeatletter

\def\endequation{\eqno \hbox{\@eqnnum}$$\global\@ignoretrue}
\def\@eqnnum{{\normalfont \normalcolor (\theCor)}}
 \makeatother

 \renewcommand{\labelenumi}{\rm (\arabic{enumi})} 

 \newcommand\rest[1]{_{{\textstyle|}#1}}

 \newcommand\w{\textstyle{\bigwedge^2}}
 \DeclareMathSymbol{\boxtimes}{\mathbin}{AMSa}{"02}
 \newcommand\wave{\widetilde}
 \newcommand\Span[1]{\left<#1\right>}
 \newcommand\Spancomma{\Span{\ ,\ }}
 \newcommand\Sum{\bigoplus}
 \newcommand\ve{{}^{\vee}}
 \newcommand\1{^{-1}}
 \newcommand\iso{\cong}
 \newcommand\broken{\mathrel{{\relbar\kern-.2pt\rightarrow}}}
 \newcommand\ot{\leftarrow}
 \newcommand\isoto{{@>{\approx}>>}}
 \newcommand{\tensor}{\otimes}
 \newcommand\into{\hookrightarrow}
 \DeclareMathSymbol{\onto}{\mathrel}{AMSa}{"10}
 \newcommand\C{\mathbb C}
 \newcommand\proj{\mathbb P}
 \newcommand\R{\mathbb R}
 \newcommand\Z{\mathbb Z}
 \newcommand\Ga{\Gamma}
 \newcommand\De{\Delta}
 \newcommand\La{\Lambda}
 \newcommand\Om{\Omega}
 \newcommand\al{\alpha}
 \newcommand\be{\beta}
 \newcommand\ga{\gamma}
 \newcommand\de{\delta}
 \newcommand\la{\lambda}
 \newcommand\fie{\varphi}
 \newcommand\si{\sigma}
 \newcommand\om{\omega}
 \newcommand\fm{{\mathfrak m}}
 \newcommand\fS{{\mathfrak S}}
 \newcommand\sE{{\mathcal E}}
 \newcommand\sF{{\mathcal F}}
 \newcommand\sK{{\mathcal K}}
 \newcommand\sL{{\mathcal L}}
 \newcommand\sM{{\mathcal M}}
 \newcommand\Oh{{\mathcal O}}
 \newcommand\sP{{\mathcal P}}
 \newcommand\sV{{\mathcal V}}
 \newcommand\Aut{\operatorname{Aut}}
 \newcommand\Coh{\operatorname{Coh}}
 \newcommand\sEnd{\operatorname{{\mathcal E\kern-0.4mm\it nd}}}
 \newcommand\Ext{\operatorname{Ext}}
 \newcommand\Grass{\operatorname{Grass}}
 \newcommand\Hom{\operatorname{Hom}}
 \renewcommand\Im{\operatorname{Im}}
 \newcommand\Jac{\operatorname{Jac}}
 \newcommand\SO{\operatorname{SO}}
 \newcommand\Pfaff{\operatorname{Pfaff}}
 \newcommand\Pic{\operatorname{Pic}}
 \newcommand\Pico{\operatorname{Pic^0}}
 \newcommand\Res{\operatorname{Res}}
 \newcommand\Supp{\operatorname{Supp}}
 \newcommand\sli{\operatorname{{\mathfrak{sl}}}}
 \newcommand\SL{\operatorname{SL}}
 \newcommand\SU{\operatorname{SU}}
 \newcommand\U{\operatorname{U}}
 \newcommand\GL{\operatorname{GL}}
 \newcommand\ho{\operatorname{hom}}
 \newcommand\im{\operatorname{im}}
 \newcommand\coker{\operatorname{coker}}
 \newcommand\rank{\operatorname{rank}}
 \newcommand\sgn{\operatorname{sgn}}
 \newcommand\ch{\mathrm{ch}}
 \newcommand\ev{\mathrm{ev}}
 \newcommand\id{\mathrm{id}}
 \newcommand\Lloc{L^{\mathrm{loc}}}
 \newcommand\Lrat{L^{\mathrm{rat}}}
 \newcommand\sLloc{{\mathcal L}^{\mathrm{loc}}}
 \newcommand\sLrat{{\mathcal L}^{\mathrm{rat}}}
 \newcommand\Mbar{\overline M}

\title{Non-Abelian Brill--Noether theory\\ and Fano 3-folds}
\author{Shigeru Mukai \thanks{Lecturer at the Sept.\ 1995 meeting of the
Math.\ Soc.\ Japan at T\^ohoku Univ., Sendai}}
\date{}

\begin{document}
\maketitle
\markboth{Brill--Noether theory and Fano 3-folds}{Shigeru Mukai}

The Jacobian variety of an algebraic curve $C$ is one connected component
of the moduli space of line bundles over $C$; the moduli space of vector
bundles can be viewed as a non-Abelian generalization of this. Here we
fix a line bundle $\xi$, and consider the moduli space of rank 2 vector
bundles with fixed determinant:
 \[
M_C(2,\xi)=\left.\left\{
\renewcommand{\arraycolsep}{2pt}
\begin{array}{l}
\text{stable rank $2$ vector}\\
\text{bundles $E$ over $C$}
\end{array}
\right | \,\w E\iso\xi
\right\}\Bigm/\text{(isomorphism)}.
 \]
The number $h^0(L)$ of linearly independent section of a line bundle $L$
can be used to define subschemes of $\Jac C$, called the {\em
Brill--Noether locuses}. These have been studied since the 19th century,
since they reflect properties of an individual curve that are beyond the
control of the Riemann--Roch theorem. In this article, we recall this
theory briefly in \S2, then generalize it to the moduli spaces
$M_C(2,\xi)$, and give applications to Fano manifolds and curves on K3s.
There are various types  of generalizations; here we treat Types~II
and~III. For Type~I, see, for example, \cite{36}, \cite{44}.
 \[
 \renewcommand{\arraystretch}{2}
 \begin{tabular}{|c|c|c|c|c|}
 \hline
 $
 \renewcommand{\arraystretch}{1}
 \begin{matrix}
 \text{algebraic}\\
 \text{group}
 \end{matrix}$
 & 
 $
 \renewcommand{\arraystretch}{1}
 \begin{matrix}
 \text{moduli space of}\\
 \text{principal bundles}
 \end{matrix}$
 & type & measure &
 $
 \renewcommand{\arraystretch}{1}
 \begin{matrix}
 \text{expected}\\
 \text{codimension}
 \end{matrix}$
 \\[5pt]
 \hline
 $\C^*$ & $\Jac C$ & classic & $h^0(L)$ & $h^0(L)h^1(L)$ \\
 \hline
 $\GL(r,\C)$ & $M_C(r,d)$ & Type I & $h^0(E)$ & $h^0(E)h^1(E)$ \\
 \hline
 $\SL(2,\C)$ & $M_C(2,\xi)$ & Type II & $\ho(F,E)$ &
$\displaystyle{\binom{\ho(F,E)}2}$\\[4pt]
 \hline
 $\SL(2,\C)$ & $M_C(2,K_C)$ & Type III & $h^0(E)$ &
$\displaystyle{\binom{h^0(E)+1}2}$\\[4pt]
 \hline
 \end{tabular}
 \]

The moduli spaces $\Jac C$ and $M_C(2,\xi)$ are K\"ahler manifolds, but
their under\-lying structures of symplectic manifolds parametrize unitary
representations of the fundamental group $\pi_1(C)$. Here $\Jac C$
parametrizes representations in $\U(1)$, and if $\xi$ has even degree then
$M_C(2,\xi)$ parametrizes irreducible representations in $\SU(2)$ (see
Narasimhan and Seshadri \cite{33}).

From now on, we let $C$ be a curve of genus $\ge2$. As is well known, the
universal cover of $C$ is the upper half-plane $H=\{\Im z>0\}$. We fix a
subgroup\footnote{This is the same thing as a theta characteristic or spin
structure on $C$, that is, the choice of an element in $\frac12K_C$.}
$\Ga\subset\SL(2,\R)$ which maps isomorphically to $\pi_1(C)$ under the
natural surjective homomorphism $\SL(2,\R)\to\Aut H$, and consider the
space $S_1(\Ga,\rho)$ of auto\-morphic forms of weight 1 with coefficients
in $\rho$. In the classic case, when $d=g-1$, the Brill--Noether locuses
$W^r_{g-1}$ are nothing other than the theta divisor $\Theta\subset\Jac C$
of the Jacobian and its singular sets; these parametrize representations
$\rho\colon\Ga\to U(1)$ for which $\dim S_1(\Ga,\rho)\ge r+1$.
 \begin{figure}[ht]
 \centering\mbox{\kern-1.5cm\epsfbox{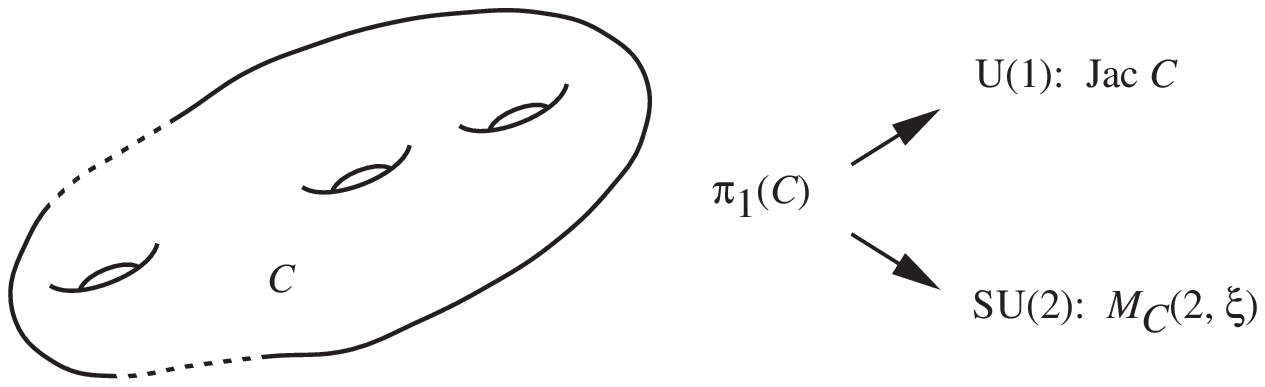}\kern-1cm}
 \end{figure}

In a similar way, the Type III Brill--Noether locuses
 \[
M_C(2,K,n)=\bigl\{ E \bigm | h^0(E)\ge n+2 \bigr\}
 \]
parametrize irreducible representations $\rho\colon\Ga\to\SU(2)$ for
which the space
 \[
S_1(\Ga,\rho)=
\left\{
 \begin{pmatrix}
 f(z)\\ g(z)
 \end{pmatrix}
\left |\,
 \begin{array}{l}
 \text{$f(z)$ and $g(z)$ are holomorphic functions}\\[2pt]
 \text{on $H$, and for all $\ga=\begin{pmatrix} a&b\\c&d
\end{pmatrix}\in\Ga$,}\\[6pt]
 \begin{pmatrix}
 f(\frac{az+b}{cz+d})\\ g(\frac{az+b}{cz+d})
 \end{pmatrix}
 =(cz+d)\rho(\ga)
 \begin{pmatrix}
 f(z)\\ g(z)
 \end{pmatrix}
 \end{array}
\right.
\right\}
 \]
of vector valued automorphic forms of weight 1 has dimension $\ge n+2$.

The study of Brill--Noether locuses in practical cases is an illustration
of the power of determinantal formalism for defining equations. \S1 recalls
the general theory of determinantal locuses. It is interesting to note that
the 3 different types of non-Abelian Brill--Noether locuses have different
types of determinantal forms. Brill--Noether locuses of Type~III are defined
by conditions on the dimension of intersection of two Lagrangian subspaces
in a symplectic vector space (see \S5). Similarly, Brill--Noether locuses
of Type~II are defined as subschemes as the set of points where a
skewsymmetric homomorphism between two dual (or twisted dual) vector
bundles drops rank (see \S6).

One of the original examples inspiring our interest in non-Abelian
Brill--Noether locuses was the linear section theorem for a general curve
of genus~8:
 \[
C=\Grass(2,6)\cap H_1\cap H_2\cap\cdots\cap H_7
 \]
(see \cite{21}, \cite{22}); to prove this, we had to construct an element
$E\in M_C(2,K)$ such that $h^0(E)\ge6$. The theorem itself has its origins
in the classification of Fano 3-folds, and Brill--Noether locuses have
many other nice applications to Fano 3-folds. The point is that certain
Fano 3-folds are realized as special cases of Brill--Noether locuses of
Type~III and Type~II. Thus for a curve of genus 7, the Brill--Noether
locus $M_C(2,K,3)=\{h^0(E)\ge5\}$ is a Fano 3-fold of genus 7 (see \S8).
Similarly, for a given stable vector bundle $F$ over a curve $C$ of genus
3, the locus $M_C(2,K\colon 3F)=\{\hom(F,E)\ge3\}$ is a Fano 3-fold of
genus 9 (see \S9). In this sense, the current article follows on from
\cite{24}.

Non-Abelian Brill--Noether locuses are also closely related to vector
bundles over K3 surfaces. It is known that any 2-dimensional component of
the moduli space of these is again a K3 surface. The search for this kind
of phenomenon over curves is a further inspiration for our study of
non-Abelian Brill--Noether locuses. In this vein, we are following on from
\cite{20}. For moduli spaces of vector bundles $E$ over a surface, we
obtain varieties of various dimensions by changing the specified value of
the second Chern class $c_2(E)$; for curves, changing the number of global
sections $h^0(E)$ or the number of homomorphisms $\hom(F,E)$ can be viewed
as a substitute for the choice of $c_2(E)$. In the final section \S10, as
an application of this idea, we try out a non-Abelian version of the
Albanese morphism; when a curve $C$ is contained in a K3 surface $S$, we
show how to recover $S$ from $C$ in terms of {\em double moduli}, that is,
a twice iterated moduli construction.

 \medskip\noindent{\bf Notation and terminology\enspace} We consider
algebraic varieties over the complex numbers $\C$. We write $V\ve$ or
$E\ve$ to denote the dual of a vector space $V$ or a vector bundle $E$, and
$\chi(E)=\sum(-1)^ih^i(E)$ for the Euler--Poincar\'e characteristic.

 \section{Determinantal subschemes}

It frequently happens that an algebraic variety can be defined by a
system of equations given as the minors of a matrix. The most basic case
of this is the degeneracy locus of a homomorphism $f\colon E\to F$ of
vector bundles over a variety $X$.

 \begin{Def} We write $D_k(f)$ for the set of zeros of the $k$th exterior
power of $f$
 \[
\bigwedge^kf \colon \bigwedge^k E \to \bigwedge^k F;
 \]
more precisely, $D_k(f)$ is the zero subscheme of $\bigwedge^kf$. We also
write $D^k(f)=D_{r-k}(f)$, where $r=\rank E$.
 \end{Def}
The homomorphism $f$ is locally expressed as a matrix with entries regular
functions on $X$, so that $D_k(f)$ is the locus (subscheme) of common
zeros of all the $k\times k$ minors of this matrix. Thus as subsets of
$X$, $D_k(f)$ and $D^k(f)$ are equal to
 \[
\bigl\{ x\in X
\bigm| \rank(f_k\colon E_x\to F_x)\le k
\bigr\}
\enspace\text{and}\enspace
\bigl\{ x\in X
\bigm| \dim \ker (f_k\colon E_x\to F_x)\ge k
\bigr\}.
 \]

 \begin{Prop}[\cite7] \label{Prop:gen_det} Suppose that $X$ is
nonsingular. Then
 \begin{enumerate}
 \item $\dim D_k(f)\ge\dim X-(r-k)(s-k)$, where $r$ and $s$ are the ranks
of $E$ and $F$; we assume here that $D_k(f)$ is nonempty.
 \item $D_{k-1}(f)$ is contained in the singular locus of $D_k(f)$.
 \item Write $\rho$ for the right hand side of the inequality of (1).
Then $D_k(f)$ is nonsingular of dimension $\rho$ at a point $x\in
D_k(f)\setminus D_{k-1}(f)$ if and only if the linear map
 \[
\ker f_x \tensor (\coker f_x)\ve \to \fm_x/\fm_x^2
 \]
corresponding to $x$ is injective.
 \item {\em The Giambelli--Thom--Porteous formula:} define the Chern
polynomial of the homomorphism $f$ by
 \[
\sum_{i\ge0} c_i(f)t^i=
\left(
\sum_{i\ge0} c_i(E)t^i
\right)\1
\left(
\sum_{i\ge0} c_i(F)t^i
\right).
 \]
where $c_i(f)\in H^{2i}(X)$. Assume that equality holds in the inequality
of (1). Then the fundamental cohomology class of $D_k(f)$ is given by the
following Schur polynomial in the Chern classes of $f$:
 \[
[D_k(f)] = \De_{\underbrace{s-k,\dots,s-k}_{r-k\ \mathrm{times}}}
\bigl(c(f)\bigr)\in H^{2N-2\rho}(X),
 \]
where $N=\dim X$.
 \end{enumerate}
 \end{Prop}

Here in (4), the Schur polynomial $\De_{\la_1,\la_2,\dots,\la_m}(c)$ is 
defined quite generally for $\la_1\ge\la_2\ge\dots\ge\la_m$ as the
determinant of the matrix having $(i,j)$th entry $c_{\la_i-i+j}$. In other
words,
 \[
\De_{\la_1,\la_2,\dots,\la_m}(c) = \det
\left|
\begin{array}{llcc}
c_{\la_1}&c_{\la_1+1}&c_{\la_1+2}&\dots\\
c_{\la_2-1}&c_{\la_2}&c_{\la_1+1} \\
\ \vdots & &\ddots\\
&&& c_{\la_m}
\end{array}
\right|.
 \]

Assume that $D_{k-1}=\emptyset$, and that the assumption in (3) holds for
all $x\in D_k(f)$. In this case, the conormal bundle to $D_k(f)$ in $X$ is
isomorphic to
 \begin{equation}
 \ker(f\rest{D(f)}) \tensor \coker(f\rest{D(f)})\ve.
 \label{eq:1.3}
 \end{equation}

In cases when the matrix representing a homomorphism of vector bundles is
symmetric or skewsymmetric, the dimension of the degeneracy locus is
bigger than the $\rho$ in Proposition~\ref{Prop:gen_det}, so that (3) and
(4) are not meaningful as they stand. However, we can modify them to make
them meaningful. Consider first the symmetric case. Perhaps the obvious
way to get symmetry is to take a homomorphism $f$ from a vector bundle $E$
to its dual $E\ve$ satisfying $f\ve=f$. In this case, the degeneracy locus
$D^k(f)$ satisfies the following.

 \begin{Prop} \label{Prop:sym_det}
 \begin{enumerate}
 \item $\dim D^k(f)\ge\dim X-k(k+1)/2$. Here we assume that $D^k(f)$ is
nonempty.
 \item $D^{k+1}(f)$ is contained in the singular locus of $D^k(f)$.
 \item Write $\si$ for the right hand side of the inequality in (1).
Then $D^k(f)$ is nonsingular of dimension $\si$ at a point $x\in
D^k(f)\setminus D^{k+1}(f)$ if and only if the linear map
 \[
S^2\ker f_x \to \fm_x/\fm_x^2
 \]
corresponding to $x$ is injective.
 \item {\em The Harris--Tu formula \cite{11}:} If equality holds in the
inequality of (1) then the fundamental class of $D^k(f)$ is given by the
Schur polynomial:
 \[
[D^k(f)]=2^k\De_{k,k-1,\dots,2,1}\bigl(c(E\ve)\bigr)\in H^{2N-k(k+1)}(X),
 \]
 \end{enumerate}
where $N=\dim X$.
 \end{Prop}

A case when the appearance of symmetry is perhaps slightly less obvious is
that of a pair of Lagrangian subbundles. As a warming up exercise,
consider a finite dimensional vector space $V$ with a skew inner product,
that is, a nondegenerate skew\-symmetric bilinear form
 \[
 \Spancomma\colon V\times V \to \C
 \]
A vector subspace $L\subset V$ with $2\dim L=\dim V$ such that the skew
inner product $\Spancomma$ restricts to zero on $L$ is called a {\em
Lagrangian subspace}. Note that the quotient space $V/L$ is then dual to
$L$ under the skew inner product $\Spancomma$. Given a pair $(L_1,L_2)$ of
Lagrangian subspaces, the composites
 \[
 L_1\into V\onto V/L_2\iso L_2\ve \quad\text{and}\quad
 L_2\into V\onto V/L_1\iso L_1\ve
 \]
provide us with two linear maps. Now clearly, these are dual to one another,
and both have the same kernel $L_1\cap L_2$.

 \begin{Ex} If $L_1$ is a vector space and $L_1\ve$ its dual, we can give
the direct sum $V=L_1\oplus L_1\ve$ the standard skew inner product
 \[
 \begin{pmatrix}
0&I_n\\ -I_n&0
 \end{pmatrix},
\quad\text{where $n=\dim L_1$.}
 \]
The two direct summands $L_1$ and $L_1\ve$ are both Lagrangian subspaces.
Now let $L_2$ be an $n$-dimensional subspace of $V$ with $L_1\ve\cap L_2=0$;
any such $L_2$ is the graph of a linear map $f\colon L_1\ve\to
L_1$,
 \[
\Ga_f=\bigl\{ (f(x),x) \bigm| x\in L_1\ve \bigr\}.
 \]
We easily check the following:
 \begin{enumerate}
 \item $L_1\cap L_2\iso\ker f$;
 \item $L_2=\Ga_f$ is a Lagrangian subspace if and only if $f$ is
symmetric.
 \end{enumerate}
 \end{Ex}

To globalize the above construction, we now consider a vector bundle $\sV$
together with a skew inner product
 \[
 \Spancomma\colon\sV\times\sV \to \Phi
 \]
with values in a line bundle $\Phi$. A Lagrangian vector subbundle
$\sL\subset\sV$ is defined as above, and we consider a pair $(\sL_1,\sL_2)$
of such.

 \begin{Def} Consider the composite of inclusion and quotient maps
 \[
\sL_1 \into \sV \onto\sV/\sL_2.
 \]
We write $D^k(\sL_1,\sL_2)$ for its degeneracy locus; as a set, it is the
locus of $x\in X$ such that $\dim\sL_{1,x}\cap\sL_{2,x}\ge k$.
 \end{Def}

 \begin{Ex}\label{Ex:1.7} Consider the direct sum $\sV$ of a vector bundle
$\sL$ and its dual $\sL\ve$, with the standard skew inner product. The
graph of a symmetric homo\-morphism $f\colon\sL\ve\to\sL$ is a Lagrangian
subbundle $\Ga_f\subset\sV$, and for every $x\in X$ we have
$\sL_x\cap\Ga_{f,x}\iso\ker f_x$. Thus $D^k(\sL,\Ga_f)$ coincides with
$D^k(f)$.
 \end{Ex}

Locally, any pair $(\sL_1,\sL_2)$ of Lagrangian vector subbundles can be
written in terms of a symmetric homomorphism as in Example~\ref{Ex:1.7}, so
that $D^k(\sL_1,\sL_2)$ is the set of common zeros of the minors of a
symmetric matrix.

 \begin{Prop} \label{Prop:lang_det}
 \begin{enumerate}
 \item $\dim D^k(\sL_1,\sL_2)\ge\dim X-k(k+1)/2$; here we assume that
$D^k(\sL_1,\sL_2)$ is nonempty.
 \item $D^{k+1}(\sL_1,\sL_2)$ is contained in the singular locus of
$D^k(\sL_1,\sL_2)$.
 \item Write $\si$ for the right hand side of the inequality of (1). Then
$D^k(\sL_1,\sL_2)$ is nonsingular of dimension $\si$ at a point $x\in
D^k(\sL_1,\sL_2)\setminus D^{k+1}(\sL_1,\sL_2)$ if and only if the linear
map
 \[
\Phi_x\1\tensor S^2(\sL_{1,x}\cap\sL_{2,x})\to\fm_x/\fm_x^2
 \]
corresponding to $x$ is injective.
 \end{enumerate}
 \end{Prop}

Formulas for the fundamental class of the degeneracy locus
$D^k(\sL_1,\sL_2)$ have been obtained by Pragacz and Fulton
(\cite{38}, \cite{39}); the former uses the Schur $Q$-polynomials.

Assume that $D^{k+1}(\sL_1,\sL_2)=\emptyset$, and that the assumption in
(3) of the Proposition holds for every $x\in D^k(\sL_1,\sL_2)$; then, as in
(\ref{eq:1.3}), we get an exact sequence
 \begin{equation}
 \begin{matrix}
 \kern-1.6cm
 0\to \Om_{D(\sL_1,\sL_2)} \to \Om_X{}\rest{D(\sL_1,\sL_2)}
 \kern-1.6cm
 \\[10pt]
 &
 \kern-1.6cm \to \Phi\1{}\rest{D(\sL_1,\sL_2)}\tensor
 S^2\bigl(\sL_1{}\rest{D(\sL_1,\sL_2)}\cap
\sL_2{}\rest{D(\sL_1,\sL_2)}\bigr) \to 0.
 \kern-1.6cm
 \end{matrix}
 \label{eq:1.9}
 \end{equation}

We finally consider the skewsymmetric case. For an $n\times n$
skewsymmetric matrix $A=(a_{ij})_{1\le i,j \le n}$, with
$a_{ij}+a_{ji}=0$, it is well known that
 \[
\det A =
\begin{cases}
0 & \text{if $n$ is odd;}\\
(\Pfaff A)^2 & \text{if $n$ is even.}
\end{cases}
 \]
Here the Pfaffian of $A$ is
 \[
 \Pfaff A=\frac{1}{2^n\cdot n!}
 \sum_{\si\in\fS_{2n}}\sgn(\si)a_{\si(1)\si(2)}\cdots
a_{\si(2n-1)\si(2n)}.
 \]

Let $f\colon E\to E\ve$ be a skewsymmetric homomorphism of vector bundles,
that is, $f\ve+f=0$. The usual way of defining the degeneracy locus of $f$
in terms of minors is unsuitable, because it introduces an unwanted
nilpotent structure. The right method is to replace the minors by the
Pfaffians of principal sub\-matrixes of even order. Note first that the
rank of $f_x\colon E_x\to E_x\ve$ is automatically even for any $x\in X$.
Write $r=\rank E$; then for an integer $\nu\ge0$ with $\nu\equiv r\mod 2$,
we define
 \begin{equation}
 \begin{array}[t]{rcl}
 P^\nu(f) & = & \bigl\{ x\in X \bigm | \rank f_x\le r-\nu \bigr\};\\[3pt]
 & = & \text{locus of common zeros of the Pfaffians of}\\
 && \text{$(r-\nu)\times(r-\nu)$ principal submatrixes}\\
 && \text{of a skewsymmetric matrix representing $f$.}
 \end{array}
 \label{eq:1.10}
 \end{equation}
In this display, the first line defines $P^\nu(f)$ as a point set; the
second line only makes sense locally, but provides a system of defining
equations for $P^\nu(f)$. Thus putting the two together defines the
subscheme structure of $P^\nu(f)$.

 \begin{Prop} \label{Prop:skew_det}
 \begin{enumerate}
 \item If $P^\nu(f)$ is nonempty then
 \[
\dim P^\nu(f) \ge \dim X-\nu(\nu-1)/2.
 \]
 \item $P^{\nu+2}(f)$ is contained in the singular locus of $P^\nu(f)$.
 \item Write $\tau$ for the right hand side of the inequality of (1).
Then $P^\nu(f)$ is nonsingular of dimension $\tau$ at a point $x\in
P^\nu(f)\setminus P^{\nu+2}(f)$ if and only if
the linear map
 \[
 \w\ker f_x \to \fm/\fm_x^2
 \]
corresponding to $x$ is injective.
 \item {\em The Harris--Tu formula \cite{11}:} If equality holds in (1)
then the fundamental class of $P^\nu(f)$ is given by
 \[
[P^\nu(f)]=\De_{\nu-1,\nu-2,\dots,2,1}\bigl(c(E\ve)\bigr)\in
H^{2N-\nu(\nu-1)}(X),
 \]
where $N=\dim X$.
 \end{enumerate}
 \end{Prop}

 \begin{Rem} For a line bundle $\Phi$, we say that a homomorphism
$f\colon E\to E\ve\tensor\Phi$ is {\em twisted skew\-symmetric} if
$f\ve\otimes 1_\Phi+f=0$; then (4) holds on replacing $c(E)$ by
$c(E\otimes\sqrt{\Phi}\1)$ (the {\em squaring principle} of Harris
and Tu, see \cite{11}).
 \end{Rem}

Assume that $P^{\nu+2}=\emptyset$ and that the assumption of (3) holds for
every $x\in P^\nu(f)$. Then, as in (\ref{eq:1.3}) and (\ref{eq:1.9}), we get
the exact sequence
 \begin{equation}
0\to\Om_{P(f)} \to \Om_X{}\rest{P(f)} \to
\Phi\1{}\rest{P(f)}\tensor\w\ker(f\rest{P(f)}) \to 0.
 \label{eq:1.13}
 \end{equation}

 \begin{Rem}\label{Rem:1.14} In (4) of Proposition~\ref{Prop:skew_det},
the right hand side is a polynomial in the Chern characters $\ch_i(E)$ of
$E$, and this formula contains only the Chern characters of odd degree.
For example, when $\nu=3$ or $4$, we get
 \begin{align*}
 \De_{2,1}\bigl(c(E\ve)\bigr)
 &\hbox{} =-\left|\begin{matrix} c_2&c_3\\ 1 & c_1 \end{matrix}\right|
 =\frac{c_1^3}{3}-2\ch_3, \\
 \De_{3,2,1}\bigl(c(E\ve)\bigr)
 &\hbox{} =\left|\begin{matrix} c_3&c_4&c_5\\
 c_1&c_2&c_3 \\ 
 0&1&c_1 \end{matrix}\right|
 =\frac{1}{45}c_1^6-\frac{1}{3}c_1^3\ch_3{}+24c_1\ch_5{}-4\ch_3^2{}.
 \end{align*}
 \end{Rem}

 \begin{Rem} The 3 types of degeneracy locus in
Propositions~\ref{Prop:gen_det}, \ref{Prop:lang_det} and
\ref{Prop:skew_det} correspond to the bounded symmetric domains of
Type~I${}_{r,s}$, Type~III (Siegel upper half-space) and Type~II.
 \end{Rem}

 \section{Brill--Noether theory}

The set $H^1(\Oh_C^*)$ of isomorphism classes of all holomorphic line
bundles over a curve $C$ has a natural complex structure. This is called
the {\em Picard variety} of $C$, denoted by $\Pic C$. The degree of a line
bundle decomposes $\Pic C$ as the disjoint union of its connected
components:
 \[
\Pic^dC=\bigl\{\xi\bigm|\deg\xi=d\bigr\}\bigm/\text{(isomorphism)}.
 \]
The short exact sequence
 \[
0\to\Z @>{2\pi\sqrt{-1}}>>\Oh_C @>{\exp}>>\Oh_C^*\to1
 \]
of sheaves over $C$ induces the long exact cohomology sequence
 \[
\cdots\to H^1(C,\Z)\to H^1(\Oh_C)\to H^1(\Oh_C^*)\to H^2(C,\Z)\to0,
 \]
and one sees from this that $\Pic^dC$ is a $g$-dimensional complex torus,
where $g$ is the genus of $C$. The classic Brill--Noether locuses for
curves are defined inside $\Pic^dC$ by setting
 \[
W^r_d=\bigl\{\xi\bigm|
\text{$\deg\xi=d$ and $h^0(\xi)\ge r+1$}\bigr\}\subset\Pic^dC.
 \]
The number $h^0(\xi)$ of linearly independent global sections of $\xi$ is
an upper\-continuous function of $\Pic^dC$, so that the $W^r_d$ are closed
subsets of $\Pic^dC$. We leave the definition of the $W^r_d$ as
subschemes until later, and list their known properties (compare
\cite{1}).

 \no\label{no:2.1}
$\dim W_d^r\ge g-(r+1)(r-d+g)$. The right hand side is called the
{\em Brill--Noether number}, and is usually denoted by $\rho$.

 \no\label{no:2.2}
$W_d^{r+1}$ is contained in the set of singular points of $W_d^r$.

 \no\label{no:2.3}
Suppose that the line bundle $\xi$ is in $W_d^r\setminus W_d^{r+1}$, in
other words, that $H^0(\xi)$ is exactly $(r+1)$-dimensional. Then
$W_d^r$ is nonsingular of dimension $\rho$ at the point $[\xi]$ if and
only if the multiplication map (called the {\em Petri map})
 \[
H^0(\xi)\tensor H^0(K_C\xi\1)\to H^0(K_C) \iso H^1(\Oh_C)\ve
 \]
is injective.

 \no\label{no:2.4}
If equality holds in the inequality of (\ref{no:2.1}) then the fundamental
class of $W_d^r$ is given by
 \[
[W_d^r]=\la(r,d,g)\Theta^{g-\rho}\in H^{2g-2\rho}(\Pic^d C).
 \]
Here $\Theta$ is the theta divisor, and we set
 \[
\la(r,d,g)=\prod_{i=0}^r \frac{i!}{(g-d+r+i)!}.
 \]

 \no\label{no:2.5} If the Brill--Noether number $\rho$ is nonnegative
then $W_d^r$ is nonempty (see \cite{13}, \cite{14}). Also, if $\rho>0$
then $W_d^r$ is connected (\cite8). Conversely, for a general curve
$C\in\sM_g$, if $\rho<0$ then $W_d^r=\emptyset$.

We give the point set $W_d^r\subset\Pic^dC$ a subscheme structure as the
degeneracy locus of a homomorphism of vector bundles. We choose distinct
points $P_1,\dots,P_N$ on $C$, and define a divisor $D$ as their sum,
$D=\sum_{i=1}^N P_i$. We choose $N$ sufficiently large so that
$H^1(\xi(D))=0$ for every $\xi\in\Pic^dC$. Now consider the short exact
sequence of sheaves on $C$
 \[
0\to\xi\to\xi(D)\to\Sum_{i=1}^N\xi(D)\rest{P_i}\to0
 \]
and the induced cohomology long exact sequence
 \begin{equation}
0\to H^0(\xi) \to H^0(\xi(D)) @>{\al_\xi}>> \Sum_{i=1}^N
H^0(\xi(D)\rest{P_i})\to H^1(\xi) \to 0.
 \end{equation}
We obviously have $h^0(\xi)=h^0(\xi(D))-\rank\al_\xi$. Now globalizing
this, we get a homomorphism of vector bundles over $\Pic^dC$ of the form
 \begin{equation}
A=\coprod_{\xi\in\Pic C}\al_\xi\colon\: E=\coprod_{\xi\in\Pic C}
H^0(\xi(D)) \longrightarrow F=\Sum_{i=1}^NF_i,
\label{eq:2.7}
 \end{equation}
where $F_i= \coprod_{\xi\in\Pic C}H^0(\xi(D)\rest{P_i})$. To say this in a
more precise way, let $\sL$ be a Poincar\'e line bundle. That is, $\sL$ is
a line bundle on the product $C\times\Pic^dC$ such that
$\sL\rest{C\times[\xi]}\iso\xi$ for every $\xi\in\Pic^dC$. Consider the
exact sequence of sheaves
 \[
0\to\sL\to\sL\tensor_{\Oh_C}\Oh_C(D)\to\sL\tensor_{\Oh_C}\Oh_D(D)\to0
 \]
on the direct product; then taking direct images with respect to the
second projection $\pi\colon C\times\Pic^dC\to\Pic^dC$ induces the
cohomology long exact
sequence
 \[
0\to\pi_*\sL\to
\pi_*(\sL\tensor_{\Oh_C}\Oh_C(D))@>A>>\pi_*(\sL\tensor_{\Oh_C}\Oh_D(D))\to
R^1\pi_*\sL\to0;
 \]
the homomorphism $A$ appearing here gives the precise definition of
(\ref{eq:2.7}). In fact, $A$ restricted to the fiber over each point
$[\xi]\in\Pic^dC$ is exactly the $\al_\xi$ of (\ref{eq:2.7}).

 \begin{Def} We define the subscheme structure of the Brill-Noether locus
$W_d^r\subset\Pic^dC$ as the degeneracy locus $D^{r+1}(A)$ of the above
homomorphism of vector bundles
 \[
A\colon
\pi_*(\sL\tensor_{\Oh_C}\Oh_C(D))\to\pi_*(\sL\tensor_{\Oh_C}\Oh_D(D)).
 \]
 \end{Def}

The assertions (\ref{no:2.1}), (\ref{no:2.2}), (\ref{no:2.3}) follow
easily from this definition and Proposition~\ref{Prop:gen_det}. (But to
get $W_d^r\ne\emptyset$ when $\rho\ge0$, we also have to prove that
$\la(r,d,g)\ne0$.) By the Grothendieck--Riemann--Roch theorem, one sees
that the Chern class polynomial $\sum_{i\ge0}c_i(A)t^i$ of the vector
bundle homomorphism $A$ is equal to $\exp(t\Theta)$; therefore a
calculation based on Proposition~\ref{Prop:gen_det}, (4) gives that the
fundamental class of
$W_d^r$ is
 \begin{align*}
[W_d^r] &\hbox{} =
\De_{\underbrace{r-d+g,\dots,r-d+g}_{\text{$r+1$ times}}}
(c_i=\Theta^i/i!)\\ &\hbox{} = \la(r,d,g)\Theta^{g-\rho}\in
H^{2g-2\rho}(\Pic^dC).
 \end{align*}
This is (\ref{no:2.4}).

 \begin{Rem} The scheme structure of $W_d^r$ just defined is independent
of the choice of the auxiliary divisor $D=\sum_{i=1}^N P_i$. In fact, by
definition, $W_d^r$ is the subscheme defined by a Fitting ideal of the
cokernel of the homomorphism $A$, which is the first higher direct
image sheaf $R^1\pi_*\sL$; however, the Fitting ideal of a module is
independent of its realization as the cokernel of a homomorphism of
vector bundles (see \cite{37}).
 \end{Rem}

 \section{The moduli space of stable vector bundles}

The set of isomorphism classes of all holomorphic vector bundles of rank
$r$ over a curve $C$ is equal to the cohomology set $H^1(\GL(r,\Oh_C))$.
However, in contrast to the case of line bundles ($r=1$), before we can
give this set a complex structure, we first have to restrict the set of
vector bundles under study.

 \begin{Def}[Mumford \cite{27}] A rank $2$ vector bundle $E$ over a
curve $C$ is {\em stable} if
 \[
\deg\xi<\frac{1}{2}\deg E
 \]
for every line subbundle $\xi\subset E$. Moreover, if the inequality
holds in the weaker form $\le$, we say that $E$ is {\em semistable}.
 \end{Def}

We consider the moduli space
 \[
M_C(2,\xi)=\left.\left\{
\renewcommand{\arraycolsep}{2pt}
\begin{array}{l}
\text{stable rank $2$ vector}\\
\text{bundles $E$ over $C$}
\end{array}
\right | \,\w E\iso\xi
\right\}\Bigm/\text{(isomorphism)}.
 \]
of stable vector bundles over $C$ with fixed determinant line bundle
$\xi$. The first thing to note is that the set of all first order
infinitesimal deformations of $E$ is parametrized by the vector space
$H^1(\sEnd E)$. Here $\sEnd E$ is the sheaf formed by the (local)
endomorphisms of $E$, or in other words, $\sEnd E=E\ve\tensor E$.
However, $\sEnd E$ splits as a direct sum
 \[
\Oh_C\cdot\id\oplus \sli E,
 \]
where the first summand consists of scalar multiplication by functions, and
the second of endomorphisms having zero trace. Simple considerations lead
to the following facts:

 \no\label{no:3.2}
Infinitesimal deformations that leave $\det E$ invariant are parametrized
by the subspace $H^1(\sli E)$. Stability of $E$ is an open condition,
preserved by small perturbations; thus $H^1(\sli E)$ is the tangent space
to the moduli space $M_C(2,\xi)$ at the point $[E]$. Obstructions to
deformations live in the space $H^2(\sli E)$, which is zero for reasons
of dimension. Also, by the Riemann--Roch theorem, $H^1(\sli E)$ has
dimension $3g-3$, so that we obtain the following result:

 \no
$M_C(2,\xi)$ is nonsingular and of dimension $3g-3$.
 \medskip

We now consider the global structure. By Geometric Invariant Theory
(\cite{29}), $M_C(2,\xi)$ is a quasiprojective algebraic variety (see
\cite{35}, \cite{42}). Up to iso\-morphism, it only depends on the parity
of $\deg\xi$. For odd degree, $M_C(2,\xi)$ is projective, and in
particular compact; its second cohomology group is isomorphic to $\Z$,
and its first Chern class $c_1$ equal to twice the positive generator. Thus
$M_C(2,\xi)$ is a Fano manifold of dimension $3g-3$ and index 2 (see
\cite{40}). For a hyperelliptic curve $C$, it has the following concrete
description.

 \begin{Thm}[Desale--Ramanan \cite{4}] Let $C$ be the hyperelliptic curve
defined by the equation $y^2=\prod_{i=1}^{2g+2}(x-\la_i)$. Then the odd
moduli space $M_C(2,\xi)$ is the closed subset of the Grassmannian
variety $\Grass(g,2g+2)$ consisting of $(g-1)$-planes of $\proj^{2g+1}$
contained in the complete intersection of two quadrics
$V_4\subset\proj^{2g+1}$ defined by
 \[
V_4:
 \left(
 \sum_{i=1}^{2g+2} X_i^2 = \sum_{i=1}^{2g+2} \la_iX_i^2 = 0
 \right)
\subset\proj^{2g+1}.
 \]
 \end{Thm}

The even moduli space $M_C(2,\xi)$, that is, when $\xi$ has even degree,
is not compact. However, it is contained as an open set in a projective
algebraic variety $\Mbar_C(2,\xi)$, with points of the boundary
parametrizing the isomorphism classes of bundles that are direct sums of
the form
 \[
L\oplus \xi L\1\quad\text{for $L\in\Pic^dC$,\quad
where $d=\deg\xi/2$.}
 \]
In more rigorous terms, $\Mbar_C(2,\xi)$ is the
moduli space of so-called $S$-equivalence classes of semistable bundles.
It is singular along the boundary (the case of genus 2 is an exception);
however, the anticanonical line bundle exists, and $\Mbar_C(2,\xi)$
is a singular Fano manifold of index 4 (see (\ref{eq:4.4}) below).

We can construct a natural map from the even moduli space $M_C(2,\xi)$ to
projective space as follows. Write $d=\deg\xi/2$, and set
 \[
D_E=\bigl\{\eta \bigm | H^0(E\tensor\eta\1)\ne0 \bigr\}
\subset\Pic^{d+1-g}C\quad\text{for $E\in M_C(2,\xi)$}.
 \]
This is a codimension 1 subvariety, and as a divisor, it is an element of
the linear system $|2\Theta|$; we thus obtain a morphism
 \[
M_C(2,\xi)\to\proj^{2^g-1}=|2\Theta|,\quad
\text{given by}
\quad E\mapsto D_E.
 \]
This morphism extends naturally to the compactification $\Mbar_C(2,\xi)$
(see \cite{30} or \cite{Beauville}).

 \begin{Exs} \begin{enumerate}
 \item For a genus $2$ curve $C$, the above morphism is an iso\-morphism
$\Mbar_C(2,\xi)\isoto\proj^3$, and the boundary of the
compactification maps to the Kummer quartic surface associated with the
Jacobian variety $\Jac C$.
 \item Consider the above morphism
 \[
\Mbar_C(2,\xi)\to\proj^7
 \]
for a genus $3$ curve $C$. If $C$ is hyperelliptic, then it is a double
covering of the quadric $Q^6$. For nonhyperelliptic $C$, that is, when $C$
is a plane quartic, it is an embedding, with image a quartic hypersurface
(see \cite{32}). This quartic hypersurface is singular along the Kummer
\hbox{$3$-fold} of $\Jac C$, and its defining equation can be derived from
this fact (see \cite2).
 \end{enumerate}
 \end{Exs}

We now explain the Hecke correspondence between the even and odd moduli
space (compare \cite{31} and \cite{Beauville}). For a rank 2 vector bundle
$E$ and a point $p$, suppose that we are given a 1-dimensional subspace of
the fiber $E_p$; write $E_p\to k(p)$ for the linear map with this as
kernel. Now consider its composite
 \[
E\onto E_p\onto k(p)
 \]
with the natural projection map $E\to E_p=E/\fm_pE$. Its kernel $F$ is a
vector bundle over $C$; we say that $F$ is obtained by {\em pinching} $E$
at $p$. We have $\w F\iso (\w E)(-p)$, so that $\deg F=\deg E-1$. Taking
the dual of the exact sequence defining $F$
 \[
0\to F\to E\to k(p)\to0
 \]
gives an exact sequence
 \[
0\to E\ve\to F\ve\to k(p)\to0.
 \]
Thus conversely, we can also recover $E\ve$ by pinching $F\ve$ at $p$. To
express this contruction geometrically, given the $\proj^1$-bundle
$\proj(E)$ and a point $x$ on it, we first blow up $\proj(E)$ at $x$,
then contract the birational transform $f'$ of the fiber $f$ through $x$
to a point $y$, thus obtaining a new $\proj^1$-bundle $\proj(F)$.
 \begin{figure}[ht]
 \centering\mbox{\kern-1cm\epsfbox{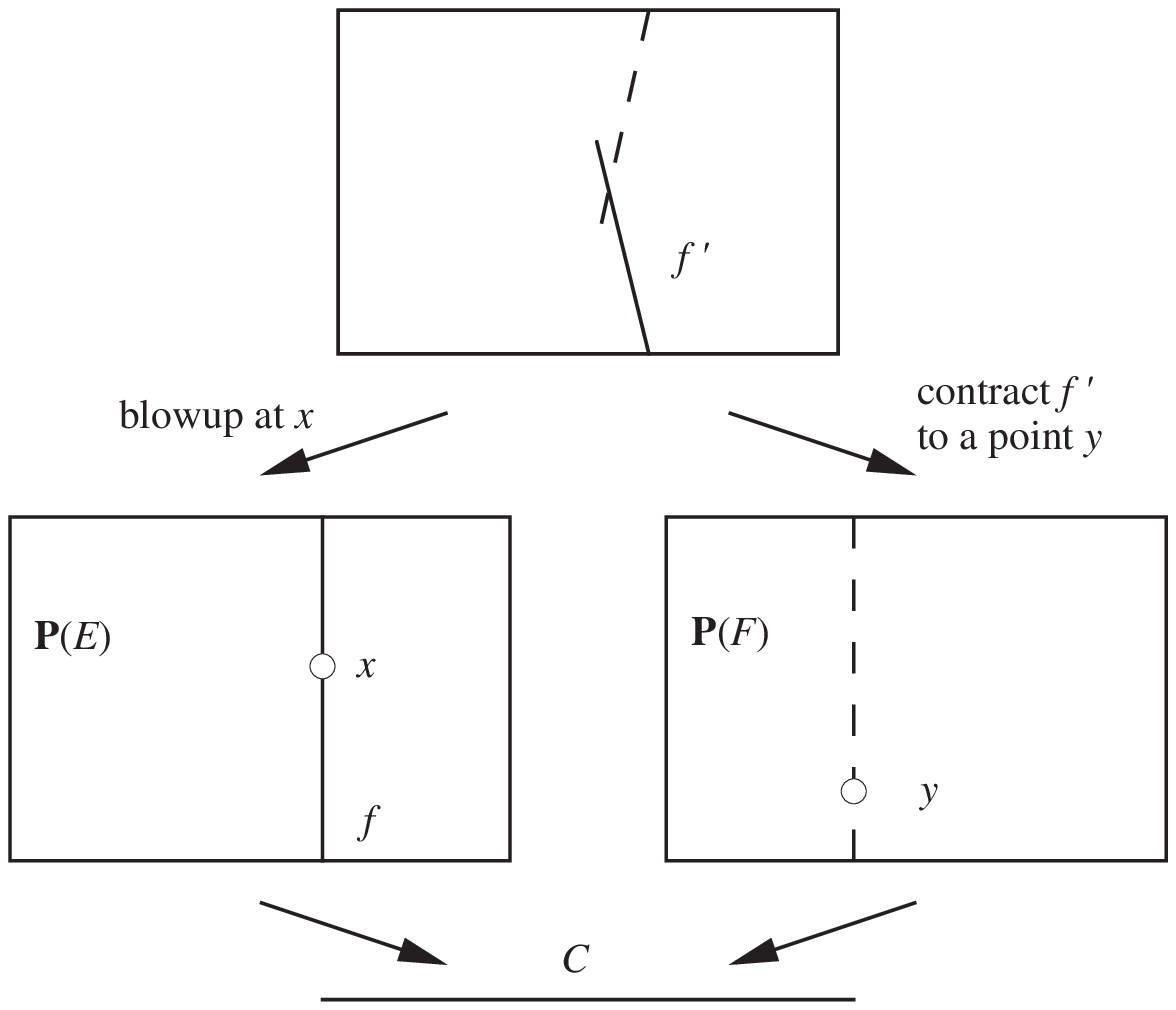}\kern-1cm}
 \end{figure}

We now fix a line bundle $\xi$ of odd degree and a point $p$ of $C$.
Moreover, let $\sE$ be the universal vector bundle over the direct product
$C\times M_C(2,\xi)$, and write $\sE_p$ for its restriction to $p\times
M_C(2,\xi)$. For each $E\in M_C(2,\xi)$ we have
$\sE\rest{C\times[E]}\iso E$, so that the fiber of $\sE_p$ at $[E]$ is
isomorphic to the fiber $E_p$ of $E$ at $p$. We write
$Z=\proj(\sE_p)$ for the $\proj^1$-bundle over $M_C(2,\xi)$ obtained by
projectivizing $\sE_p$; it parametrizes pairs consisting of a stable
vector bundle $E\in M_C(2,\xi)$ together with a 1-dimensional vector
subspace of the fiber $E_p$ at $p$. For a point $Z\ni z=(E,f_z\colon
E_p\to\C)$, we write $F_z$ for the vector bundle obtained by pinching $E$
at $p$ as described above. $F_z$ is defined by the exact sequence
 \[
0\to F_z\to E@>{\wave f_z}>> k(p)\to0.
 \]
Now since $E$ is stable and of odd degree, $F_z$ is semistable. We
therefore obtain a classifying map
 \[
\fie\colon Z\to \Mbar_C(2,\xi(-p))
 \]
from $Z$ to the moduli space, defined by $z\mapsto[F_z]$. The
correspondence
 \begin{equation}
 \renewcommand{\arraystretch}{1.4}
 \begin{matrix}
 Z&@>\fie>>&\Mbar_C(2,\xi(-p))\\
 \kern-2cm\text{$\proj^1$-bundle}
 \downarrow
 \hphantom{\text{$\proj^1$-bundle}}\kern-2cm\\
 \kern-1cm M_C(2,\xi) \kern-1cm
 \end{matrix}
 \label{eq:3.6}
 \end{equation}
obtained in this way is called the {\em Hecke correspondence}. In what
follows, we refer to $Z$ as the {\em Hecke graph}; it can alternatively
be viewed as a moduli space of stable parabolic bundles (compare \cite{41},
\cite{15}).

 \no
Over the open subset $M_C(2,\xi(-p))$, every fiber of $\fie$ is isomorphic
to the projective line $\proj^1$. (Over boundary points, the behaviour of
$\fie$ provides interesting examples of extremal rays.)

 \no\label{no:3.8}
There exists a vector bundle $\sF$ over the direct product $C\times Z$ 
with $\sF\rest{C\times z}\iso F_z$ for each point $z\in Z$. This is
obtained from the universal bundle $\sE$ over $C\times M_C(2,\xi(-p))$ by
taking the pullback $(1\times\pi)^*\sE$ and then pinching it as follows:
because $Z=\proj(\sE_p)$, it has a tautological line bundle $\Oh_Z(H)$,
together with a tautological surjective sheaf homomorphism
$\pi^*\sE_p\to\Oh_Z(H)\to0$. Noting that $\sE_p=\sE\rest{p\times M}$, we
need only take $\sF$ to be the kernel of the composite homomorphism
 \[
(1\times\pi)^*\sE\to\pi^*\sE\rest{Z\times p}\to\Oh_Z(H).
 \]

 \section{Non-Abelian Brill--Noether locuses\\ of Type~III}

The even moduli space admits two particular cases, when $\xi$ is equal to
the trivial bundle $\Oh_C$ or the canonical line bundle $K_C$.

 \begin{Thm}[Narasimhan--Seshadri \cite{33}, Donaldson \cite6]
For a rank $2$ vector bundle $E$ over $C$, the following $2$ conditions
are equivalent.

 \begin{enumerate}
 \item $E\in M_C(2,\Oh_C)$, or in other words, $E$ is stable and $\w
E\iso\Oh_C$.
 \item $E$ is a flat bundle obtained from an irreducible $\SU(2)$
representation of the fundamental group $\pi_1(C)$.
 \end{enumerate}

In particular, $M_C(2,\Oh_C)$ can be identified with the set of
equivalence classes of irreducible representations of $\pi_1$ in $\SU(2)$.
 \end{Thm}

The other particular case, the moduli space $M_C(2,K)$ of bundles with
the canonical bundle as determinant, parametrizes representations
$\rho\colon\Ga\to\SU(2)$, where $\Ga\subset\SL(2,\R)$ is a lift of $\pi_1$
(see the introduction).

 \begin{Prop} The space $H^0(E)$ of global sections of a bundle $E\in
M_C(2,K)$ is isomorphic to the space $S_1(\Ga,\rho_E)$ of automorphic
forms with coefficients in the corresponding representation $\rho_E$.
 \end{Prop}

We use the dimension of this vector space to define the Brill--Noether
locuses:
 \begin{equation}
M_C(2,K,n)= \bigl\{ [E] \bigm | h^0(E)\ge n+2 \bigr\} \subset M_C(2,K).
 \end{equation}
In the case $n=-1$, we write $\Xi$ for $M_C(2,K,-1)$. This is a
codimension 1 subvariety, and we call the corresponding line bundle
$\Oh_M(\Xi)$ the {\em determinant line bundle} of $M_C(2,K)$.
Its fiber at the point $[E]$ is naturally isomorphic to the 1-dimensional
vector space
 \[
\det H^0(E)\tensor \det H^1(E)\ve.
 \]
Moreover, the canonical divisor of $M_C(2,K)$ is given by
 \begin{equation}
 K_M\sim -4\Xi,
 \label{eq:4.4}
 \end{equation}
and this extends as a linear equivalence of Cartier divisors to the
compactification $\Mbar_C(2,K)$.

We give $M_C(2,K,n)$ a subscheme structure in essentially the same way as
for $W_d^r$. We fix an effective divisor $D$ of degree $\ge g-1$; then for
every $E\in M_C(2,K)$, we have
 \[
H^1(E(D))=H^0(E(-D))=0.
 \]
The short exact sequence of sheaves on $C$
 \begin{equation}
0 \to E \to E(D) \to E(D)\rest{D} \to 0
 \label{eq:4.5}
 \end{equation}
induces a cohomology long exact sequence
 \begin{equation}
0 \to H^0(E) \to H^0(E(D)) @>{\al_E}>> H^0(E(D)\rest{D}) \to H^1(E) \to 0,
 \label{eq:4.6}
 \end{equation}
so that $h^0(E)=2\deg D-\rank\al_E$. Next, we globalize $\al_E$ to a
homomorphism
 \[
A\colon \coprod_E H^0(E(D))\to \coprod_E H^0(E(D)\rest{D})
 \]
of vector bundles over the moduli space $M_C(2,K)$, and define the
subscheme structure of $M_C(2,K,n)$ as its degeneracy locus $D^{n+2}(A)$.
(We omit the technical question of how to deal with the fact that the
universal bundle does not exist.)

Applying Proposition~\ref{Prop:gen_det} to $D^{n+2}(A)$ gives the
estimate: if $M_C(2,K,n)$ is nonempty, then
 \[
\dim M_C(2,K,n)\ge 3g-3-(n+2)^2.
 \]
However, this is pretty feeble; the right estimate is the following:

 \begin{Thm}\label{Thm:4.7}
 \begin{enumerate}
 \item If $M_C(2,K,n)$ is nonempty then
 \[
\dim M_C(2,K,n)\ge 3g-3-(n+2)(n+3)/2.
 \]

 \item Write $\si$ for the right hand side of this inequality. Moreover,
suppose that for some $E\in M_C(2,K)$, the space $H^0(E)$ has dimension
exactly $n+2$. Then $M_C(2,K,n)$ is nonsingular at $[E]$ of dimension
$\si$ if and only if the linear map defined by multiplication (which we
call the {\em Petri map})
 \[
S^2H^0(E) \to H^0(S^2E) \iso H^1(\sli E)\ve
 \]
is injective
(see (\ref{no:3.2}) for $H^1(\sli E)$).
 \end{enumerate}
 \end{Thm}

To prove the theorem, we need to consider $M_C(2,K,n)$ with a different
structure. We explain this in the following section. The analog of the
classic existence result (\ref{no:2.5}) for $M_C(2,K,n)$ is an unsolved
problem:

 \begin{Prob} 
 \begin{enumerate}
 \item Suppose that in Theorem~\ref{Thm:4.7}, (1), the right hand side
is $\ge0$; then is it true that $\Mbar_C(2,K,n)\ne\emptyset$? (Here we
put $\Mbar$ to allow semistable bundles.)
 \item Does equality holds in (1) of the theorem for a general curve
$C\in\sM_g$?
 \end{enumerate}
 \end{Prob}

The cases in which $M_C(2,K,n)$ is well understood are shown in the
following table (compare also \cite{Laszlo}, \cite{PR}, \cite{OPP}). Here
the curve $C$ is assumed to satisfy suitable generality assumptions, which
we omit. Case~(6) of the table is described in more detail in \S8.
 \[
 \renewcommand{\arraystretch}{2}
 \begin{tabular}{|c|c|c|c|c|l|}
 \hline
 & $g(C)$ &
 $
 \renewcommand{\arraystretch}{1}
 \begin{matrix}
 \dim M_C(2,K)\\=3g-3
 \end{matrix}
 $
 & n+2 & dimension &
 $
 \renewcommand{\arraystretch}{1}
 \renewcommand{\arraycolsep}{0pt}
 \begin{array}{l}
 M_C(2,K,n)\\
 \text{and its boundary}
 \end{array}
 $
 \\
 \hline
(1) & 3 & 6 & 2 & 3 & 
 $
 \renewcommand{\arraystretch}{0.9}
 \renewcommand{\arraycolsep}{0pt}
 \begin{array}{l}
 \text{cone over Veronese $\proj^2$,}\\
 \text{the boundary is $\iso C$}
 \end{array}
 $
 \\
 \hline
(2) & 3 & 6 & 3 & 0 & single point \\
 \hline
(3) & 4 & 9 & 3 & 3 & singular cubic 3-fold \\
 \hline
(4) & 5 & 12 & 4 & 2 & $\proj^2\supset$ quintic curve \\
 \hline
(5) & 6 & 15 & 5 & 0 & single point \\
 \hline
(6) & 7 & 18 & 5 & 3 & Fano 3-fold of genus 7 \\
 \hline
(7) & 8 & 21 & 6 & 0 & single point \\
 \hline
(8) & 9 & 24 & 6 & 3 & singular quartic 3-fold \\
 \hline
(9) & 11 & 30 & 7 & 2 & polarized K3 of genus 11 \\
 \hline
 \end{tabular}
 \]

 \section{Lagrangian formalism for Brill--Noether \\ locuses}

Let $D$ be an effective divisor of degree $\ge g-1$. Instead of the exact
sequence (\ref{eq:4.5}), we consider
 \begin{equation}
0 \to E(-D) \to E(D) \to A_E \to 0.
 \end{equation}
The cokernel $A_E$ is an Artinian sheaf with the same support as $D$, and
$V_E=H^0(A_E)$ is a vector space of dimension $4\deg D$. Consider the
skewsymmetric bihomomorphism
 \begin{equation}
 E(D) \times E(D) \to \w E(D)\iso K_C(2D),
 \label{eq:5.2}
 \end{equation}
and the induced map
 \begin{equation}
 A_E \times A_E \to K(2D)/K.
 \end{equation}
We take global sections, and compose with the residue map
 \[
 \renewcommand{\arraystretch}{1.3}
 \begin{matrix}
 H^0(K(2D)/K) & \to & H^1(K_C)\iso \C \\
 (\om_p)_p & \mapsto & \sum_{p\in\Supp D} \Res_p \om_p
 \end{matrix}
 \]
to obtain a skewsymmetric bilinear map
 \begin{equation}
 \renewcommand{\arraystretch}{1.3}
 \begin{matrix}
\Spancomma\colon& V_E\times V_E & \to & H^1(\w E)\iso \C, \\[3pt]
 &\bigl((f_p)_p,(g_p)_p\bigr) & \mapsto
 & \sum_{p\in\Supp D} \Res_p f_p\wedge g_p.
 \end{matrix}
 \end{equation}
Now the original pairing (\ref{eq:5.2}) was nondegenerate, so that
$\Spancomma$ is also nondegenerate. Thus $V_E$ is a symplectic vector space;
moreover, it contains two naturally occurring Lagrangian subspaces. One of
these is
 \begin{equation}
 \Lloc:= H^0(E/E(-D)),
 \end{equation}
which is obviously Lagrangian by definition. The other is
 \begin{equation}
 \Lrat:= \im \bigl[ H^0(E(D)) \to V_E = H^0(E(D)/E(-D)) \bigr].
 \end{equation}
First, the Riemann--Roch formula gives
 \begin{align*}
 \dim\Lrat &\hbox{} = h^0(E(D)-h^0(E(-D))\\
 &\hbox{} = \chi(E(D))=2\deg D,
 \end{align*}
that is, exactly half the dimension of $V_E$. The fact that the
skewsymmetric pairing $\Spancomma$ restricts to zero on $\Lrat$ is a
consequence of the following standard result.

 \begin{RThm} If $\om$ is a rational (Abelian) differential on a curve
$C$ then the sum of its residues is zero: $\sum_p\Res_p\om=0$.
 \end{RThm}

The following statement is obvious:

 \begin{Prop}[\cite{23}] The image of $H^0(E)\to V_E$ is exactly the
intersection of the two Lagrangian subspaces described above:
 \[
\im\bigl[H^0(E)\to V_E\bigr]=\Lloc\cap\Lrat.
 \]
 \end{Prop}

Now $\Lloc\cap\Lrat$ is the kernel of the composite homomorphism
 \[
\Lrat\into V \to V/\Lloc \iso (\Lloc)\ve,
 \]
but we note that this map is nothing other than $\al_E$ in (\ref{eq:4.6}).

Globalizing the above argument is routine. We assume for the sake of
convenience that the universal bundle exists. Then we can construct a
vector bundle $\sV$ over the moduli space $M_C(2,K)$, together with a
skewsymmetric bihomo\-morphism $\Spancomma\colon\sV\times\sV\to\Phi$ and two
Lagrangian subbundles $\sLloc,\sLrat\subset\sV$, so that at each point
$[E]\in M_C(2,K)$, the fiber is $V_E$, together with $\Lloc,\Lrat$ as
above. Taking $\deg D$ to be sufficiently large, we get the following:

 \begin{Prop} At each point $[E]$ of the moduli space $M_C(2,K)$, we have
 \[
H^0(E)\iso\sLloc_{[E]}\cap\sLrat_{[E]};
 \]
therefore the Brill--Noether locus $M_C(2,K,n)$ coincides with the
degeneracy locus $D^{n+2}(\sLloc,\sLrat)$ of the pair $(\sLloc,\sLrat)$
of Lagrangian subbundles of $\sV$.
 \end{Prop}

Thus Theorem~\ref{Thm:4.7} follows from Proposition~\ref{Prop:lang_det}.
Moreover, from the exact sequence (\ref{eq:1.9}) and a calculation, we
get the following result.

 \begin{Prop} \label{Prop:5.10} If $n$ is odd then
$M_C(2,K,n+1)=\emptyset$. Assume in addition that the Petri map
$S^2H^0(E)\to H^0(S^2E)$ is injective for every $E\in M_C(2,K,n)$. Then
the canonical bundle of
$M_C(2,K,n)$ is isomorphic to the restriction of $\Oh_M((n-5)/2)$; here
$\Oh_M(1)$ is the determinant line bundle of the moduli space $M_C(2,K)$.
 \end{Prop}

\section{Non-Abelian Brill--Noether locuses of Type~II}

We say that a pair $(F,E)$ of rank 2 vector bundles over a curve $C$ has
{\em canonical difference} if $\w F\tensor K_C\iso\w E$. We want to
consider Brill--Noether locuses for such pairs, defined in terms of the
space of homomorphisms $\Hom(F,E)$. In this vein, we could of course allow
both of $F$ and $E$ to vary, but here we work with a fixed $F$. To start
with, we write the dual bundle $F\ve$ any-old-how as an extension of two
line bundles, that is, we fix an exact sequence
 \begin{equation}
0\to M\to F\ve\to N\to0.
 \label{eq:6.1}
 \end{equation}
Now tensoring this with $E$ and taking global sections gives a long exact
sequence
 \begin{equation}
 \begin{aligned}
 0\to\hbox{}&H^0(M\tensor E)\to\Hom(F,E)\to H^0(N\tensor E)\\
 @>\de_E>>\hbox{}&H^1(M\tensor E)\to\Ext^1(F,E)\to H^1(N\tensor E)\to0.
 \end{aligned}
 \label{eq:6.2}
 \end{equation}
The natural bihomomorphism
 \[
(N\tensor E)\times(M\tensor E)\to N\tensor M\tensor\w E\iso K_C
 \]
induces a cup product
 \[
\Spancomma\colon H^0(N\tensor E)\times H^1(M\tensor E)\to H^1(K_C)\iso\C,
 \]
which is nondegenerate by Serre duality. On the other hand, a simple
calculation with \v Cech cocycles shows that
 \[
\Span{s,\de_E(s)}=0 \quad\text{for all $s\in H^0(N\tensor E)$}.
\]
Thus the connecting homomorphism $\de_E$ in (\ref{eq:6.2}) is
skewsymmetric. In particular, $\rank\de_E$ is even, so that we obtain the
following.

 \begin{Prop} Let $(F,E)$ be a pair of rank $2$ vector bundles over a
curve $C$ with canonical difference; then
 \[
\dim\Hom(F,E)\equiv\deg F\mod2.
 \]
 \end{Prop}

Now, for any integer $\nu\ge0$ with $\nu\equiv\deg F\mod2$, we define a
subset of $M_C(2,\w F\tensor K_C)$ by
 \[
M_C(2,K\colon \nu F):= \bigl\{ E \bigm | \dim\Hom(F,E)\ge\nu \bigr\}.
 \]
We can put a subscheme structure on $M_C(2,K\colon \nu F)$ by using the
globalized form of (\ref{eq:6.2}). For convenience, we explain this under
the pretence that there is a universal bundle $\sE$ over $C\times
M_C(2,\w F\tensor K_C)$, although of course no such thing exists when
$\deg F$ is even. Then pulling back the exact sequence (\ref{eq:6.1}) to
the direct product and tensoring with $\sE$ gives a sequence
 \[
0 \to \sE\tensor_{\Oh_C}M
 \to \sE\tensor_{\Oh_C}F\ve
 \to \sE\tensor_{\Oh_C}N
 \to 0
 \]
over $C\times M_C(2,\w F\tensor K_C)$. Taking direct images by the
second projection $\pi$ gives
 \begin{equation}
 \begin{aligned}
0 \to\hbox{}& \pi_*(\sE\tensor_{\Oh_C}M)
 \to \pi_*(\sE\tensor_{\Oh_C}F\ve)
 \to \pi_*(\sE\tensor_{\Oh_C}N)\\
 @>{\De}>>\hbox{}& R^1\pi_*(\sE\tensor_{\Oh_C}M)
 \to R^1\pi_*(\sE\tensor_{\Oh_C}F\ve)
 \to R^1\pi_*(\sE\tensor_{\Oh_C}N) \to 0.
 \end{aligned}
 \label{eq:6.4}
 \end{equation}
The connecting homomorphism $\De$ restricted to every fiber equals the
$\de_E$ of (\ref{eq:6.2}), which is skewsymmetric with respect to Serre
duality. Twisting by a suitable line bundle $\Phi$, we can make $\De$
itself skewsymmetric with respect to Serre duality (relative to $\pi$):
 \[
R^1\pi_*(\sE\tensor_{\Oh_C}M)\iso\pi_*(\sE\tensor_{\Oh_C}N)\ve\tensor\Phi;
 \]
compare the start of \S7 for the twist $\Phi$.

 \begin{Def} 
We write $M_C(2,K\colon \nu F)\subset M_C(2,\w F\tensor K_C)$ for the subscheme 
defined as the Pfaffian locus $P^\nu(\De)$ (see (\ref{eq:1.10})) of the
twisted skew\-symmetric vector bundle homomorphism
$\De\colon \pi_*(\sE\tensor_{\Oh_C}N) \to R^1\pi_*(\sE\tensor_{\Oh_C}M)$, and
call it a non-Abelian Brill--Noether locus of Type~II.
 \end{Def}

 \begin{Rem} To prove that this definition is independent of the choice
of the extension (\ref{eq:6.1}) seems not to be all that easy, but it
can be proved using a Lagrangian formalism for inner product bundles.
 \end{Rem}

By applying Proposition~\ref{Prop:skew_det} to the above definition, we get
the following result.

 \begin{Thm}
 \begin{enumerate}
 \item If $M_C(2,K\colon \nu F)\ne\emptyset$ then
 \[
\dim M_C(2,K\colon \nu F)\ge 3g-3 -\nu(\nu-1)/2.
 \]
 \item Write $\tau$ for the right hand side of this inequality. Assume also
that $E\in M_C(2,\w F\tensor K_C)$ is such that $\Hom(F,E)$ is exactly
$\nu$-dimensional.  Then $M_C(2,K\colon \nu F)$ is nonsingular of dimension
$\tau$ at $[E]$ if and only if the linear map (called the {\em Petri map})
 \[
\w \Hom(F,E) \to \Hom(\w F,S^2E)
 \]
is injective.
 \end{enumerate}
 \end{Thm}

Moreover, from the exact sequence (\ref{eq:1.13}), we get the following
result.

 \begin{Prop} \label{Prop:6.8}
Suppose that $\nu$ is odd, and $M_C(2,K\colon (\nu+2)F)=\emptyset$; suppose
also that the Petri map is injective for all $E\in M_C(2,K\colon\nu F)$.
Then the canonical class of $M_C(2,K\colon\nu F)$ is isomorphic to the
restriction of $\Oh_M((\nu-5)/2)$. Here $\Oh_M(1)$ is the positive
generator of the Picard group of the moduli space.
 \end{Prop}

 \section{Computing the fundamental class of a \\ non-Abelian
Brill--Noether locus} We want to generalize to $M_C(2,K,n)$ and
$M_C(2,K\colon\nu F)$ the formula (\ref{no:2.4}) for the fundamental class
of the locus $W^r_d$ of special line bundles. For this, we explain some
facts we need concerning the cohomology of the moduli space $M_C(2,\xi)$.
What follows is taken from Newstead \cite{34}.

First consider the odd moduli space $M_C(2,\xi)$, which is compact, and
for which the universal bundle $\sE$ exists. For each point $[E]\in
M_C(2,\xi)$ we have $\w\sE\rest{C\times[E]}\iso\w E\iso\xi$, so that
$\w\sE=\Phi\boxtimes\xi$ for some line bundle $\Phi$. Write $\phi$ for the
first Chern class of $\Phi$. Write out the Chern classes of $\sE$
separated into their K\"unneth components:
 \begin{align*}
 c_1(\sE) = \hbox{}& \phi+df\in H^2(M)\oplus H^2(C); \\
 c_2(\sE) = \hbox{}& \chi + \psi + \om\tensor f
 \end{align*}
for some $\chi\in H^4(M)$, $\psi\in H^3(M)\tensor H^1(C)$ and $\om\in
H^2(M)$, where $d=\deg\xi$ and $f\in H^2(C)$ is the fundamental class.
The only ambiguity in the choice of the universal bundle $\sE$ is
tensoring with a line bundle on $M_C(2,\xi)$, so that the 2 cohomology
classes
 \begin{align*}
 \al = \hbox{}& 2\om-d\phi\in H^2(M); \\
 \be = \hbox{}& \phi^2-4\chi\in H^4(M)
 \end{align*}
are independent of the choice of $\sE$. Also, $\psi$ is not itself a
cohomology class of $M_C(2,\xi)$, but setting
 \[
\psi^2=\ga\tensor f
 \]
gives a well-defined cohomology class $\ga\in H^6(M)$. For our purposes in
what follows, the 3 classes $\al,\be,\ga$ are all we need. The moduli space
has complex dimension $3g-3$, so that the intersection number
$(\al^m\cdot\be^n\cdot\ga^p)$ is defined when $m+2n+3p=3g-3$. These
numbers have been computed explicitly, to give the following result.

 \begin{Thm}[Thaddeus \cite{43}, Zagier \cite{48}] \label{Thm:7.1}
 \[
\bigl(\al^m\cdot\be^n\cdot\ga^p\bigr)=
(-1)^{p-g}\frac{g!m!}{(g-p)!q!} 2^{2g-2-p}(2^q-2)B_q,
 \]
where $q=m+p+1-g$ and $B_q$ is the $q$th Bernoulli number, defined by
 \[
\frac{x}{e^x-1}=1-\frac12x+\sum_{q\ \mathrm{even}} B_q\frac{x^q}{q!}.
 \]
We set $B_0=1$ and $B_q=0$ for $q<0$.
 \end{Thm}

The cohomology class $\al$ is the first Chern class of the positive
generator of the Picard group of $M_C(2,\xi)$; thus when $n=p=0$, the
number
 \[
\bigl(\al^{3g-3}\bigr)=
(-1)^g\frac{(3g-3)!}{(2g-2)!} 2^{2g-2}(2^{2g-2}-2)B_{2g-2}
 \]
is the degree of the odd moduli space.

We now compute the fundamental class of a Brill--Noether locus of Type~II,
under the assumption that it has the expected codimension. The whole point
is the connecting homomorphism (\ref{eq:6.4})
 \[
 \De\colon \pi_*(\sE\tensor_{\Oh_C}N) \to R^1\pi_*(\sE\tensor_{\Oh_C}M).
 \]
This vector bundle homomorphism is skewsymmetric under
a nondegenerate pairing
 \[
\pi_*(\sE\tensor_{\Oh_C}N) \times R^1\pi_*(\sE\tensor_{\Oh_C}M) \to \Phi
 \]
with values in the line bundle $\Phi$ defined above. We can determine the
Chern characters of the direct image by applying the
Grothendieck--Riemann--Roch formula to the virtual vector bundle
$\sE\tensor_{\Oh_M}\sqrt{\Phi}\1$, obtaining
 \begin{equation}
 \begin{array}{l}
 \ch_{2n-1}\bigl(\pi_*(\sE\tensor_{\Oh_C}N)\tensor\sqrt{\Phi}\1\bigr)
 = \\[10pt]
 \kern2cm \left\{
\begin{array}{ll}
-\al/2 & \text{for $n=1$,} \\[3pt]
\displaystyle{\left(\frac{\be}{4}\right)^{n-2}
\left(-\frac{\al\be}{8}+\frac{n-1}{2}\ga\right)\Bigm/(2n-1)!}
 & \text{for $n\ge2$.}
\end{array}
\right.
 \end{array}
 \label{eq:7.2}
 \end{equation}

The evenly numbered Chern characters can also be calculated, but are not
required because of Remark~\ref{Rem:1.14}. By
Proposition~\ref{Prop:skew_det}, (4) we get the following result.

 \begin{Prop} \label{Prop:7.3}
Suppose that the Type~II Brill--Noether locus $M_C(2,K\colon 3F)$
has codimension $3$ in $M_C(2,K_C\w F)$. Then its fundamental
class equals
 \[
(\al^3-\al\be+4\ga)/24.
 \]
 \end{Prop}

 \begin{Ex} \label{Ex:7.4}
 If the genus of $C$ is $2,3,4,5$, then under the assumption of
Proposition~\ref{Prop:7.3}, $M_C(2,K\colon 3F)$ has degree (relative to $\al$)
equal to $1$, $16$, $2544$, $1231616$.
 \end{Ex}

Next, we consider the even moduli space $\Mbar_C(2,\xi(-p))$. 
Rather than the Brill--Noether locus inside this, we take its pullback to
the Hecke graph $Z$ of (\ref{eq:3.6}), and determine the fundamental
class of this pullback. By definition, we have $c_1(\sE_p)=\al$ and
$c_2(\sE_p)=(\al^2-\be)/4$. Thus writing $H$ for the tautological divisor
of $Z=\proj(\sE_p)$, we find that $H^2=\al H-(\al^2-\be)/4$. (We continue
to write $\al,\be,\ga$ for their pullbacks to $Z$.) Using this, together
with Theorem~\ref{Thm:7.1} and the recurrence formula for the Bernoulli
numbers, we get the following result.

 \begin{Prop} \label{Prop:7.5} For $m+2n+3p=3g-3$, the intersection number
in the Hecke graph $Z$ is given by
 \[
\bigl(H\cdot\al^m\cdot\be^n\cdot\ga^p\bigr)=
\begin{cases}
\displaystyle{(-1)^{p+g}2^g\frac{g!m!}{(g-p)!q!}B_q}
 & \text{if $q\ne0$,}\\[10pt]
\displaystyle{(-1)^{p+g+1}2^g\frac{g!m!m}{(g-p)!}} & \text{if $q=0$;}
\end{cases}
 \]
here $q$ and $B_q$ are as in Theorem~\ref{Thm:7.1}.
 \end{Prop}

By (\ref{no:3.8}), we get
 \begin{align*}
c_1(\sF) & \hbox{} = \al+(d-1)f\in H^2(C\times Z),\\
c_2(\sF) & \hbox{} = \chi+\psi+ \bigl(H+(d-1)\al/2\bigr)\tensor f\in
H^4(C\times Z),
 \end{align*}
for the Chern classes of the universal vector bundle $\sF$ over $C\times
Z$. As in (\ref{eq:7.2}), we get
 \begin{equation}
 \begin{array}{l}
 \ch_{2n-1}\bigl(\pi_*(\sF(-\al/2)\tensor_{\Oh_C}N)\bigr) = \\[10pt]
 \kern2cm \left\{
\begin{array}{ll}
-H & \text{for $n=1$,} \\[6pt]
\displaystyle{\left(\frac\be4\right)^{n-2}
\left(-\frac{\be H}{4}+\frac{n-1}{2}\ga\right)\Bigm/(2n-1)!}
 & \text{for $n\ge2$.}
\end{array}
\right.
 \end{array}
 \end{equation}
(Thus in effect, we just substitute $H$ for $\al/2$ in (\ref{eq:7.2}).)
Similarly to Proposition~\ref{Prop:7.3}, we get the following result.

 \begin{Prop} \label{Prop:7.7}
 Suppose that the Type~II Brill--Noether locus $M_C(2,K\colon 4F)$
has codimension $6$ in $M_C(2,K_C\w F)$. Then its pullback to $Z$
has fundamental class equal to
 \[
\frac{h^6}{45}+\frac{h^3\ga}{18}-\frac{\ga^2+h^4\be}{36}
-\frac{h\be\ga}{45}+\frac{h^2\be^2}{180}.
 \]
 \end{Prop}

The pullback to $Z$ of the linear system $|\al|$ on $M_C(2,\xi)$ defines a
map from $Z$ to $\Mbar_C(2,\xi(-p))$ which is a double cover over the
generic point. Thus by Proposition~\ref{Prop:7.5}, one can calculate the
degree of the Brill--Noether locus with respect to the determinant line
bundle.

 \begin{Ex}
 If the genus of $C$ is $4,5,6$, then under the assumption of the
Proposition~\ref{Prop:7.7}, the degree of $M_C(2,K\colon 4F)$ equals $6$, $256$, $28640$.
 \end{Ex}

The fundamental class of Brill--Noether locus of Type~III can be
determined by applying Pragacz's formula \cite{38} to their Lagrangian
representation (see \S5).

 \begin{Prop} 
For $n=0,1$ or $2$, suppose that $M_C(2,K,n)$ has codimension exactly
$\binom{n+3}2$ in the moduli space $M_C(2,K)$. Then its pullback to
the Hecke graph $Z$ has fundamental class given by
 \begin{align*}
 & (1)\quad \frac{h^3}{6}-\frac{h\be}{6}+\frac\ga3
 \quad\text{when $n+2=2$;} \\[2pt]
 & (2)\quad \frac{h^6}{360}-\frac{h^4\be}{72}+\frac{h^3\ga}{36}
+\frac{h^2\be^2}{90}-\frac{2h\be\ga}{45}-\frac{\ga^2}{18}
 \quad\text{when $n+2=3$;} \\[2pt]
 & (3)\quad \frac{h^{10}}{302400}-\frac{h^8\be}{20160}+\frac{h^7\ga}{10080}
-\frac{17h^4\be^3}{60480}+\frac{17h^3\be^2\ga}{10080}
+\frac{h^2\be\ga^2}{720}\\[4pt]
 &\qquad +\frac{h^2\be^4}{8400} +\frac{h\ga^3}{216}-\frac{h\be^3\ga}{1050}
-\frac{3\be^2\ga^2}{2800}+\frac{h^6\be^2}{4800}-\frac{h^5\be\ga}{1200}
 \quad\text{when $n+2=4$.}
 \end{align*}
 \end{Prop}

 \section{Fano 3-folds of genus 7} The Brill--Noether locus
$M_C(2,K,3)\subset M_C(2,K)$, that is, the set of $E$ such that
$h^0(E)\ge5$, has expected codimension 15. By Proposition~\ref{Prop:5.10},
if the codimension is as expected, the anticanonical line bundle is the
restriction of the determinant line bundle of $M_C(2,K)$. In particular,
it is a Fano variety. We consider the case when $C$ has genus~7.

 \begin{Thm} \label{Thm:8.1} Let $C$ be a curve of genus $7$ not having a
special line bundle $g^1_4$. (This just means that $C$ cannot be expressed
as a cover of the Riemann sphere of degree $\le4$.) Then the Type~III
Brill--Noether locus $M_C(2,K,3)$ is a Fano $3$-fold of Picard number $1$
and genus $7$.
 \end{Thm}

We study this variety from the following two points of view:

 \begin{enumerate}
 \item the construction of birational maps and the intermediate Jacobian
variety;
 \item the degeneration of Fano 3-folds and moduli spaces.
 \end{enumerate}

There are a number of ways of constructing birational maps; here we
explain the simplest, which is obtained from 2 points $p,q\in C$. Note
that for an element $[E]\in M_C(2,K,3)$, we have $h^0(E(-p-q))\ge5-4=1$,
that is, $E$ has a nonzero global section with zeros at $p$ and $q$.
If we assume that this section has no other zeros, we obtain an exact
sequence
 \begin{equation}
 0\to\Oh_C(p+q)\to E\to K_C(-p-q)\to0.
 \label{eq:8.2}
 \end{equation}
Conversely, consider the question of when a vector bundle $E$ obtained
as an extension of the line bundles $K_C(-p-q)$ and $\Oh_C(p+q)$ has
$h^0(E)\ge5$. By Serre duality, the vector space of extensions
$\Ext^1(K_C(-p-q),\Oh_C(p+q))$ is dual to $H^0(K_C^2(-2p-2q))$. Thus if we
write
 \[
\Phi\colon C\into\proj^{13}=\proj^*H^0(K_C^2(-2p-2q))
 \]
for the embedding defined by $K_C^2(-2p-2q))$, a nonzero extension class
$e\in\Ext^1(K_C(-p-q),\Oh_C(p+q))$ determines a point of this $\proj^{13}$.
We write
 \[
\de_e\colon H^0(K_C(-p-q))\to H^1(\Oh_C(p+q)).
 \]
for the connecting homomorphism in the cohomology sequence of the
exact sequence (\ref{eq:8.2}) determined by $e$. The condition
$h^0(E)\ge5$ is equivalent to saying that $\de_e$ has rank $\le1$.
Moreover, the linear map
 \[
\Ext^1(K_C(-p-q),\Oh_C(p+q)) \to \Hom(H^0(K_C-p-q),H^1(\Oh_C(p+q)))\\
 \]
given by $e\mapsto\de_e$ is the Serre dual of the multiplication map
 \[
\mu\colon H^0(K_C(-p-q))\tensor H^0(K_C(-p-q))\to H^0(K^2(-2p-2q)).
 \]
Now consider the diagram
 \[
 \renewcommand{\arraystretch}{1.3}
 \begin{matrix}
 C & @>\quad\psi\quad>> & \proj^{13}\\
\kern-2cm\Phi_{K_C(-p-q)}
\downarrow
\hphantom{\Phi_{K_C(-p-q)}}\kern-2cm
&&
\kern-1cm\hphantom{\proj^*\mu}\downarrow{\proj^*\mu}\kern-1cm
\\
\proj^4& @>>\text{Veronese map}> & \proj^{14}
 \end{matrix}
 \]
The condition $\rank\de_e=1$ means that the point $\proj^*\mu(p)$ is
in the image of the Veronese map. Thus the set of nontrivial extensions
$E$ with $h^0(E)\ge5$ is parametrized by the quadric hypersurface
$Q^3\subset\proj^4$ which passes through the image
$C_{10}=\Phi_{K_C(-p-q)}(C)$. The vector bundle
$E$ correponding to a general point of this quadric $Q^3$ is stable, and
we obtain a birational map
 \[
M_C(2,K,3) \broken Q^3.
 \]
Taking a closer look at this correspondence, we find that this birational
map factors as follows.

 \no 
Write $J_{p,q}$ for the set of $E\in M_C(2,K,3)$ with $h^0(E(-p-q))\ge2$;
it is a rational curve of degree 2 with respect to the anticanonical class
(or the determinant line bundle). $J_{p,q}$ is clearly contained in the
indeterminacy locus of the birational map $M_C(2,K,3)\broken Q^3$, so that
we blow up $M_C(2,K,3)$ along it, writing $\wave M_C(2,K,3)$ for the
blowup.

 \no 
A vector bundle $E$ obtained as the extension corresponding to a point $x$
of the curve $C_{10}\subset Q^3\subset\proj^4$ contains the line subbundle
$K_C(-p-q-x)$, and is not stable. Therefore, it is contained in the
indeterminacy locus of the inverse birational map $Q^3\broken M_C(2,K,3)$.
Thus we write $\wave Q^3$ for the variety obtained by blowing up the
quadric hypersurface $Q^3$ along $C_{10}$.

 \no \label{no:8.5}
The two blowups $\wave M_C(2,K,3)$ and $\wave Q^3$ are isomorphic
outside a subset of codimension 2. In fact, $\wave M_C(2,K,3)$ is
obtained from $\wave Q^3$ by flopping the birational transforms
of the lines contained in $Q^3$ which are trisecant lines of $C_{10}$
(there are 14 of these in general):
 \[
 \kern2cm
 \begin{matrix}
\wave M_C(2,K,3) & \ot & \overset{\text{flops}}\cdots & \to & \wave Q^3\\
 & \searrow && \swarrow \\
\kern-1cm\text{blowup}\downarrow\hphantom{\text{blowup}}\kern-1cm&&
X_6=(2)\cap(3)\subset\proj^5 &&
\kern-1cm\hphantom{\text{blowup}}\downarrow\text{blowup}\kern-1cm \\[10pt]
\kern-3cm
M_C(2,K)\supset M_C(2,K,3)\supset J_{p,q}
\kern-2cm &&&&
\kern-2cm
C_{10}\subset Q^3\subset \proj^4
\kern-2cm &&&&
 \end{matrix}
 \]
Here $X_6$ is the common anticanonical model of $\wave M_C(2,K,3)$ and
$\wave Q^3$, so that these are two different small resolutions of the
singular Fano 3-fold $X_6$.

Because flops do not change the intermediate Jacobian, we get the
following result.

 \begin{Thm} \label{Thm:8.6}
 $M_C(2,K,3)$ is rational, and its intermediate Jacobian is isomorphic
as a principally polarized Abelian variety to the Jacobian of $C$.
 \end{Thm}

In particular, the Torelli theorem holds for $M_C(2,K,3)$. Moreover,
it can be shown that every nonsingular Fano 3-fold of genus 7 and
Picard number 1 is isomorphic to a $M_C(2,K,3)$. By the above theorem,
their moduli space is isomorphic to the moduli space of curves of
genus 7 not having a $g^1_4$.

Although we restricted Theorem~\ref{Thm:8.1} to the nonsingular case, 
even when singular, the moduli space $M_C(2,K,3)$ extends in a flat
family, provided it remains 3-dimensional. Moreover, its singularities
are all Gorenstein. There are two main directions of degenerations:
 \begin{enumerate}
 \renewcommand{\labelenumi}{(\alph{enumi})}
 \item when $C$ has an ordinary double point;
 \item when $C$ remains nonsingular, but acquires a $g^1_4$, that is,
it specializes to a 4-fold cover of the Riemann sphere.
 \end{enumerate}
In either case, we obtain a Fano 3-fold of genus 7 having an ordinary double
point. In the two cases, the (algebraic) local ring of the double point is
a unique factorization domain in Case~(a), and not so in (b).

We take a closer look at Case~(a). A curve $C$ of genus 7 with an ordinary
double point is obtained by identifying two points $p,q$ on a nonsingular
curve $\wave C$ of genus 6. In the moduli space $M_{\wave C}(2,K(p+q))$
of stable rank~2 bundles over $\wave C$, consider the locus
$M_{\wave C}(2,K(p+q),3)$ of bundles having $h^0(E)\ge5$.

 \begin{Prop} Suppose that the genus $6$ curve $\wave C$ is not trigonal,
and write $\wave C_8\subset\proj^3$ for the image of $\wave C$ under
$\Phi_{K(-p-q)}$. Then $M_C(2,K(p+q),3)$ is obtained by blowing up
$\proj^3$ in $\wave C_8$, then flopping the birational transforms of its
$5$ quadrisecants.
 \end{Prop}

The curve $\wave C_8\subset\proj^3$ is contained in a cubic surface.
After the flops, its birational transform becomes isomorphic to
$\proj^1\times\proj^1$, and can be contracted to an ordinary double
point. The singular manifold so obtained is the Brill--Noether locus
of the curve in (a). However, the node itself corresponds not to a
vector bundle over $C$, but to a torsion free sheaf.

 \section{Fano 3-folds of genus 9}

A nonhyperelliptic curve of genus 3 is isomorphic to a plane quartic curve.

 \begin{Thm} 
 Let $C$ be a nonsingular plane quartic curve, and $F$ a rank $2$ vector
bundle over $C$ of odd degree such that any section $S$ of the associated
$\proj^1$ bundle $\proj(F)$ has selfintersection number $S^2\ge3$. Then
the Brill--Noether locus of Type~II $M_C(2,K\colon 3F)$ is a nonsingular
Fano $3$-fold of index $1$ and genus $9$. Here the fact that $M$ is Fano
follows from Proposition~\ref{Prop:6.8}, and genus~$9$ from
Example~\ref{Ex:7.4}.
 \end{Thm}

 \begin{Rem} The vector bundle $F$ is stable; in fact $F$ is stable if and
only if every section $S$ has selfintersection number $S^2\ge1$.
 \end{Rem}

If we twist $F$ by a line bundle to get $F\tensor L$, the isomorphism
class of $M_C(2,K\colon 3F)$ does not change. Thus the Fano 3-folds
$M_C(2,K\colon 3F)$ depend only on the ruled surface $\proj(F)$.

As for the Fano 3-folds of genus 7 in the preceding section, the various
properties of a Fano 3-fold of genus 9 can be studied via vector bundles
in terms of moduli theory. We first construct a birational map. Let $S$
be a section of $\proj(F)$ with selfintersection $S^2=3$; in terms of
the exact sequence
 \begin{equation} 
0\to\xi\to F\to\eta\to0
 \label{eq:9.3}
 \end{equation}
associated with $S$, this means that $\deg\eta-\deg\xi=3$.

 \begin{Prop} 
 The locus of $E\in M_C(2,K\colon 3F)$ such that $H^0(\eta\1 E)$ is a rational
curve of degree $1$ with respect to the anticanonical divisor.
 \end{Prop}

We call this locus the line associated with $S$, and denote it by $l_S$.
Next, write $C_7\subset\proj^4$ for the image of the embedding
$\Phi\colon C\to\proj^*H^0(K\xi\1\eta)$. By Serre duality,
$\Ext^1(\eta,\xi)\iso H^0(K\xi\1\eta)\ve$, so that the class of the
extension (\ref{eq:9.3}) determines a point $p\in \proj^4$. Since $F$ is
stable, $p$ is not on $C_7$. Furthermore, since $\proj(F)$ does not have
a section of selfintersection 1, it is not contained on the secant
hypersurface of $C_7$. We write $C_S\subset\proj^3$ for the nonsingular
degree 7 space curve obtained by projecting $C_7$ from $p$.

 \begin{Thm} 
 Let $C,F$ and $S$ be as above. Write $\wave M_C(2,K\colon 3F)$ for the
Brill--Noether locus $M_C(2,K\colon 3F)$ with the line $l_S$ blown up and
$\wave\proj^3$ for $\proj^3$ with the curve $C_S$ blown up. Then these are
isomorphic outside a set of codimension~$2$, and are obtained from one
another by flops:
 {\rm \[
\kern1cm
 \begin{matrix}
\wave M_C(2,K\colon 3F) & \ot & \overset{\text{flops}}\cdots & \to
 & \wave\proj^3\\
 & \searrow && \swarrow \\
\kern-1cm\text{blowup}\downarrow\hphantom{\text{blowup}}\kern-1cm&&
X_{12}\subset\proj^8 &&
\kern-1cm\hphantom{\text{blowup}}\downarrow\text{blowup}\kern-1cm \\[10pt]
\kern-1cm
l_S \subset M_C(2,K\colon 3F) \subset M_C(2,K\La)
\kern-2cm &&&&
\kern-2cm
 \kern-1cm\hphantom {C_S\subset}\proj^3\supset C_S \kern-1cm
\kern-2cm &&&&
 \end{matrix}
 \]}
Here $X_{12}$ is the common anticanonical model, as in (\ref{no:8.5}).
 \end{Thm}

A point $x\in\proj^3$ of the projective space containing the degree 7 curve
$C_S$ determines a 3-dimensional vector subspace $U_x\subset
H^0(K\xi\1\eta)$. If $x\notin C_S$ then the evaluation homomorphism
 \[
 \ev_x\colon U_x\tensor\Oh_C\to K\xi\1\eta
 \]
is surjective. Now we consider the kernel of
 \[
\ev_x\tensor\xi\1\colon U_x\tensor\xi\1 \to K\xi^{-2}\eta,
 \]
writing $E_x$ for its dual. We always have $\dim\Hom(F,E_x)\ge3$.
Moreover, if $x$ is not on any quadrisecant of $C$ then $E_x$ is
stable. This defines a rational map
 \begin{equation}
 \proj^3 \broken M_C(2,K\colon 3F)\quad\text{given by $x\mapsto [E_x]$.}
 \end{equation}
The inverse can also be constructed, proving the rationality of
$M_C(2,K\colon 3F)$, but we omit the details.

Similarly to Theorem~\ref{Thm:8.6}, we obtain the following result.

 \begin{Cor} 
The intermediate Jacobian of $M_C(2,K\colon 3F)$ is isomorphic to the
Jacobian of the curve $C$.
 \end{Cor}

There are 3 main ways in which $M_C(2,K\colon 3F)$ degenerates. The first
two of these are similar to the degenerations of $M_C(2,K,3)$ in the
preceding section.
 \begin{enumerate}
 \renewcommand{\labelenumi}{(\alph{enumi})}
 \item when $C$ has an ordinary double point;
 \item when $C$ specializes to a hyperelliptic curve;
 \item finally, when the bundle $F$ is no longer general: namely, suppose
$F$ is stable, but its projectivization $\proj(F)$ has a section with
selfintersection 1. 
 \end{enumerate}

In any of the 3 cases, $M_C(2,K\colon 3F)$ acquires a single ordinary
double point. Its local ring is a unique factorization domain
in Case~(a) only. We discuss Case~(c) briefly. In this case, $F$
can be written as an extension
 \[
0\to\Oh_C\to F\to\La\to0,
 \]
where $\La$ is a line bundle of degree 1. Suppose that $h^0(F)=1$; then by
Serre duality, $h^0(K_CF\ve)=4$. Also, $K_C^2\La\1$ embeds $C$ into $\proj^4$.
Write $C_7\subset\proj^4$ for its image. The Grassmann map determined by the
sections of $K_CF\ve$ determines a quadric hypersurface $Q^3$ containing
$C_7$.

 \begin{Prop}
In Case (c), $M_C(2,K\colon 3F)$ is isomorphic to the anticanonical model
of the blowup of $C_7$ in $Q^3$.
 \end{Prop}

$C_7\subset Q^3$ has a single trisecant line. The birational transform of this
can be contracted to the ordinary double point of the anticanonical model,
and $M_C(2,K\colon 3F)$ is nonsingular outside this point.

 \section{The non-Abelian Albanese map\\ of a K3 surface}
We give an application of Brill--Noether locuses of Type~III to K3 surfaces,
and to the relation between curves and K3s. Harris and Mumford have computed
the Kodaira dimension of the moduli space of curves, and obtained the
following result:

 \begin{Thm}[\cite{10}] For sufficiently large $g$, a general curve
$C\in\sM_g$ is not isomorphic to a hyperplane section of a surface
$S\subset\proj^N$.
 \end{Thm}

Here we consider the case when $S$ is a K3.

 \begin{Prob}
Determine whether a given curve $C$ can be embedded in a K3 surface; if so,
analyze the set of all such embeddings.
 \end{Prob}

 \begin{Def}
If some embedding of $C$ into a K3 exists, we say that $C$ is a {\em K3 curve}
or a {\em K3 section}.
 \end{Def}

If a curve $C$ of genus $g\ge2$ is contained in a K3 surface $S$, then as
a divisor, it is nef and big on $S$. In fact $|C|$ is base point free, and
$C^2=2g-2>0$. Write $h$ for the linear equivalence class of the image of
$C\into S$. In what follows we consider the case that $h$ is {\em primitive}
in $\Pic S$, that is, not divisible by a natural number $\ge2$.

A K3 surface $(S,h)$ together with a (primitive) nef and big divisor is
called a {\em polarized K3}, and the natural number $g=h^2/2+1$ its {\em
genus}. We write $\sF_g$ for the moduli space of polarized K3 surfaces
$(S,h)$ of genus $g$. It is a \hbox{19-dimensional} quasiprojective
algebraic variety. We also write $\sP_g$ for the moduli space of pairs
$(S,C)$ consisting of a K3 surface together with a genus $g$ curve on it;
this has natural projection maps to $\sF_g$ and $\sM_g$, which we write
$\Psi_g$ and $\Phi_g$:
 \[
\begin{matrix}
&& \kern-.5cm\sP_g=\{(S,C)\}/\text{iso}\kern-.5cm && \into &&
 \coprod_{(S,h)\in\sF_g}|h|\\ [10pt]
& \kern-1cm\fie_g\swarrow\hphantom{\fie_g}\kern-1cm
&& \kern-1cm\hphantom{\Psi_g}\searrow\Psi_g\kern-1cm
&& \kern-2cm\hphantom{\text{$\proj^g$-bundle}}%
\swarrow\text{$\proj^g$-bundle}\kern-2cm \\ [10pt]
\begin{matrix}
\sM_g\supset\sK_g\\
\text{K3 curves}
\end{matrix}
&&&&
\kern-.3cm
\begin{matrix}
\sF_g=\{(S,h)\}/\text{iso}\\
\text{polarized K3s}
\end{matrix}
\kern-.3cm
\end{matrix}
 \]
By definition, the image of $\Phi_g$ is the locus of K3 curves; we write
$\sK_g$ for it, and $\fie_g$ for the restriction of $\Phi_g$ over it. The
moduli space $\sP_g$ is contained as an open subset in the $\proj^g$-bundle
over $\sF_g$ having fibers the complete linear systems $|h|\iso\proj^g$.
Therefore, its dimension is $19+g$. In particular, we get the following:
 \begin{equation}
\dim\sP_g\le\dim\sM_g\ \iff \ 19+g\le3g-3\ \iff\ g\ge11.
\label{eq:10.4}
 \end{equation}

 \begin{Con} 
For a curve of odd genus $g\ge11$, we can construct the inverse of $\fie_g$
as a rational map. \rm See (\ref{no:10.14})--(\ref{no:10.18}) for details;
for $g=11$, see also \cite{26}.
 \end{Con}

 \begin{Cor} \label{Cor:10.6}
For odd $g\ge11$, the locus of K3 curves $\sK_g$ is birationally equivalent
to a $\proj^g$-bundle over $\sF_g$.
 \end{Cor}

 \begin{Rem} \begin{enumerate}
 \item In spite of what inequality (\ref{eq:10.4}) may suggest, the
projection $\fie_{12}\colon\sP_{12}\to\sK_{12}\subset\sM_{12}$ is not even
generically finite. The reason behind this phenomenon is the existence of
the Fano 3-folds $V_{22}\subset\proj^{13}$ of genus $12$, discovered by
Iskovskikh. In fact, the general element of $\sK_{12}$ can be written as a
double linear section $C=V_{22}\cap H_1\cap H_2$, so that $C$ is contained
in infinitely many K3s, namely the pencil $V_{22}\cap(a_1H_1+a_2H_2)$ for
$(a_1:a_2)\in\proj^1$.

More generally, the values of $g$ for which $\fie_g$ is not birational are
exactly the values for which there exists a Fano 3-fold $V$ with
$\Pic V=\Z(-K_V)$ (explicitly, all odd numbers $g\le9$ and all even numbers
$g\le12$).

 \item It was proved in \cite{16} that $\fie_g$ is generically finite for
odd $g\ge11$. For genus 11, Corollary~\ref{Cor:10.6} also follows from the
theorem of Ciliberto and Miranda \cite3 on the irreducibility of the fiber
of $\fie_{11}$.
 \item We do not treat this here, but there are theorems of J.~Wahl
(\cite{45}, \cite{46}), saying that the Gauss map
$\w H^0(K_C)\to H^0(K_C^3)$ is not surjective for K3 curves $C$.
 \end{enumerate}
 \end{Rem}

The construction of the inverse, that is, recovering a K3 surface from a
K3 section $C\in\sK_g$ is a non-Abelian version of the Albanese map. We
first recall the theory of the dual of an Abelian variety. If $A$ is an
Abelian variety, the dual Abelian variety $\widehat A$ is the identity
component $\Pico A$ of the moduli space $\Pic A$ of line bundles on $A$;
there is a universal Poincar\'e line bundle $\sP$ over the direct product
$A\times\widehat A$. By definition,
$\widehat A\owns\al\mapsto\sP\rest{A\times\al}\in\Pico A$ is an isomorphism.
We can normalize $\sP$ by the condition that $\sP\rest{0\times\widehat A}$ is
trivial, then consider its restriction to the fibers in the other direction
$a\times \widehat A$.

 \begin{DThm}[\cite{28}, \S13] The correspondence
$a\mapsto\sP\rest{a\times\widehat A}$ defines an isomorphism
$A\isoto\skew5\widehat{\widehat A}$.
 \end{DThm}

 \begin{Cor} If $X$ is a nonsingular projective algebraic variety, the dual
Abelian variety of the Picard variety $\Pico X$ is the Albanese variety of
$X$. That is, there exists a morphism $\al\colon X\to(\Pico X)\widehat\ $,
which has the universal mapping property for maps from $X$ to an Abelian
variety.
 \end{Cor}

 \begin{Rem}
 \begin{enumerate}
  \item The Albanese variety is constructed analytically by picking a
basis $\om_1,\dots,\om_n\in H^0(X,\Om^1_X)$ for regular differential
1-forms on $X$, and then taking the integral
 \[
 \begin{matrix}
 X&\to&\C^n/\mathrm{(periods)}\\[6pt]
 x & \mapsto & \displaystyle{\left(\int^x\om_1,\dots,\int^x\om_n \right)}
 \end{matrix}
 \]
  \item In positive characteristic, if $\Pico X$ is nonreduced, we have to
replace $\Pico X$ by its reduced Abelian variety before taking the dual.
 \end{enumerate}
 \end{Rem}
 
We define a non-Abelian version of the above duality theorem for K3 surfaces
(see \cite{dual K3s} for details). Let $(S,h)$ be a polarized K3 surface of
odd genus $g=2n+1$. Consider the moduli space
 \[
M_S(2,h,s):= \left\{ E\, \left|
\renewcommand{\arraycolsep}{2pt}
\begin{array}{l}
\text{$E$ is a rank 2 stable sheaf, with}\\
\text{$c_1(E)=h$ and $\chi(E)=s+2$}
\end{array}
\right\}\Bigm/\text{(isomorphism)}\right..
 \]
This is nonsingular and of dimension $2(g-2s)$ (see \cite{18}, \cite{19},
\cite{20}). In particular, when $s=n=(g-1)/2$, it is a surface; we denote
it by $\widehat S=M_S(2,h,n)$. We assume in what follows that $n$ is also
odd, so that $g\equiv3\mod4$.

 \begin{Prop} Assume that $\widehat S$ is compact. (This is equivalent to
saying that all semistable sheaves are stable.)
 \begin{enumerate}
 \item $\widehat S$ is a K3 surface.
 \item There exists a polarization $\widehat h$ of $\widehat S$ having the
same degree as $h$.
 \item There exists a universal vector bundle $\sE$ over the direct product
$S\times\widehat S$ with $c_1(\sE)=\pi^*_Sh+\pi^*_{\widehat S}{\widehat h}$;
(we call this the {\em normalized universal vector bundle).}
 \end{enumerate}
 \end{Prop}

By definition, the correspondence $\al\mapsto\sE\rest{S\times\al}$ defines
an isomorphism $\widehat S\isoto\allowbreak M_S(2,h,n)$; restricting to the
fibers in the other direction gives the following.

 \begin{DThm} Let $(S,h)\in\sF_g$ be a general element, with $g\equiv3\mod4$,
and let $\sE$ be the normalized universal vector bundle. Then the
correspondence $s\mapsto\sE\rest{s\times\widehat S}$ defines an isomorphism
from $S$ to $\skew3\widehat{\widehat S}$.
 \end{DThm}

 \begin{Rem} \begin{enumerate}
 \item The correspondence $(S,h)\mapsto(\widehat S,\widehat h)$ extends to
an automorphism of the moduli space $\sF_g$. When $g=3$ (that is, for space
quartics), it is the identity map, for $g\ge7$, a map of order 2.

 \item In general, we have $\widehat S\not\iso S$ for $g\ge7$, but the
derived categories of coherent sheaves on $S$ and $\widehat S$ are
equivalent under an integral transformation. That is, as in \cite{17},
there is an equivalence of categories $\mathbf D(\Coh S)\isoto
\mathbf D(\Coh \widehat S)$.
 \end{enumerate}
 \end{Rem}

The construction of the inverse of $\fie_g\colon\sP_g\to\sK_g$ is as
follows. We assume that $C\in\sK_g$ is a general K3 curve of genus
$g\ge11$ with $g\equiv3\mod 4$.

 \no \label{no:10.14} Consider the moduli space $M_C(2,K)$ of vector
bundles over $C$, and the Brill--Noether locus of Type~III $M_C(2,K,n)$,
where $n=(g-1)/2$. Then $M_C(2,K,n)$ is a K3 surface, and the restriction
to it of the determinant line bundle of $M_C(2,K)$ is a polarization of
genus $g$. We write $(T,h_{\det})$ for this pair.

 \no 
There exists a universal vector bundle $\sE$ over the direct product
$C\times T$ such that $c_1(\sE)=\pi_C^*K_C+\pi_S^*h_{\det}$.

 \no 
For every $p\in C$, the restriction of $\sE$ to $p\times T$ is a stable
bundle, and has $\chi=n+2$.

 \no 
The correpondence $p\mapsto \sE\rest{p\times T}$ defines a classifying map
$C\to\widehat T=M_T(2,h_{\det},n)$, which is an embedding from $C$ to the
K3 surface $\widehat T$.

 \no \label{no:10.18}
Every K3 surface containing $C$ is isomorphic to $\widehat T$.
 \medskip

To summarize what we have said above, there is a {\em double moduli}
construction, starting from a moduli space of vector bundles over $C$,
giving a non-Abelian Albanese map $C\to\widehat T$; namely, $\widehat T$
is the dual K3 surface of a Brill--Noether locus of $C$.

When $g\equiv1\mod4$, we find in a similar way that $T=M_C(2,K,n)$ is a K3
surface. But now $h_{\det}$ is divisible by 2, and the universal vector
bundle $\sE$ does not exist globally over it. However, there does exist a
universal $\proj^1$-bundle corresponding to $\proj(\sE)$, and the K3
surface containing $C$ can be recovered as a moduli space of suitable
$\proj^1$-bundles over $T$ (we can equally well say principal
$\SO(3)$-bundles).

\medskip\noindent
Shigeru Mukai, \\
Graduate School of Polymathematics, Nagoya University, \\
Chikusa-ku, Furo-cho, Nagoya 464--01, Japan

\medskip\noindent
Translated by M. Reid \\
From S\=ugaku {\bf49}:1 (1997), 1--24 \\
To appear in the AMS series S\=ugaku Expositions, 1997 \\
The copyright on this translation will belong to the AMS.

\begin{thebibliography}{99}

 \bibitem{1} E. Arbarello, M. Cornalba, P. A. Griffiths and J. Harris,
Geometry of algebraic curves I, Springer, 1985

 \bibitem{2} A. Coble, Algebraic geometry and theta functions, Amer. Math. 
Soc. Colloq. Publ., {\bf10}, AMS, Providence, 1929

 \bibitem{3} C. Ciliberto and R. Miranda, On the Gaussian map for curves
of low genus, Duke Math. J. {\bf61} (1990), 417--443

 \bibitem{4} U. V. Desale and S. Ramanan, Classification of vector bundles
of rank 2 on hyperelliptic curves, Invent. Math. {\bf38} (1976),
161--185

 \bibitem{5} C. De Concini, and P. Pragacz, On the class of
Brill--Noether loci for Prym varieties, Math. Ann. {\bf302} (1995),
687--697

 \bibitem{6} S. K. Donaldson, A new proof of a theorem of Narasimhan and
Seshadri, J. Diff. Geom. {\bf18} (1983), 269--277

 \bibitem{7} W. Fulton, Intersection theory, Springer, 1984

 \bibitem{8} W. Fulton and R. Lazarsfeld, On the connectedness of
degeneracy loci and special divisors, Acta Math. {\bf146} (1981),
271--283

 \bibitem{9} D. Gieseker, Stable curves and special divisors: Petri's
conjecture, Invent. Math. {\bf66} (1982), 251--275

 \bibitem{10} J. Harris and D. Mumford, On the Kodaira dimension of the
moduli space of curves, Invent. Math. {\bf67} (1982), 23--86

 \bibitem{11} J. Harris and T. Lu, On symmetric and skewsymmetric
determinantal varieties, Topology {\bf23} (1984), 71--84

 \bibitem{12} J. Harris and T. Lu, Chern numbers of kernel and cokernel
bundles, Invent. Math. {\bf75} (1984), 467--475

 \bibitem{13} G. Kempf, Schubert methods with applications to algebraic
curves, Publ. Math. Centrum, Amsterdam, 1972

 \bibitem{14} S. Kleiman and D. Laksov, On the existence of special divisors,
Amer. J. Math. {\bf94} (1972), 431--436

 \bibitem{15} V. B. Mehta and C. S. Seshadri, Moduli of vector bundles on
curves with parabolic structures, Math. Ann. {\bf248} (1980), 205--239

 \bibitem{16} S. Mori and S. Mukai, The uniruledness of the moduli space
of curves of genus 11, LNM {\bf1016} (1983), Springer, pp.~334--353

 \bibitem{17} S. Mukai, Duality between $D(X)$ and $D\widehat X$ with its
application to Picard sheaves, Nagoya Math. J. {\bf81} (1981), 153--175

 \bibitem{18} S. Mukai, Symplectic structure of the moduli space of sheaves
on an Abelian or K3 surface, Invent. Math. {\bf77} (1984), 101--116

 \bibitem{19} S. Mukai, On the moduli space of bundles on K3 surfaces. I,
in Vector bundles on algebraic varieties (Bombay, 1984), Tata Inst. of
Fundamental Research Studies {\bf11}, Oxford Univ. Press, 1987,
pp.~341--413

 \bibitem{20} S. Mukai, Vector bundles over K3 surfaces and symplectic
manifolds, S\=ugaku {\bf32} (1987), 216--325; English translation,
S\=ugaku Exposition {\bf1} (1988), 139--174, A.M.S., Providence

 \bibitem{21} S. Mukai, Curves and symmetric spaces, Proc. Japan Acad.
{\bf68} (1992), 7--10

 \bibitem{22} S. Mukai, Curves and Grassmannians, in Algebraic geometry
and related topics, J.-H.~Yang, Y.~Namikawa and K.~Ueno (Eds),
International Press, Boston, 1993, pp.~19--40

 \bibitem{23} S. Mukai, Curves and Brill--Noether theory, in Current topics
in algebraic geometry, Math. Sci. Res. Inst. Publ. {\bf28}, Cambridge Univ.
Press, 1995, pp.~145--158

 \bibitem{24} S. Mukai, New developments in Fano threefolds --- On vector
bundles and moduli problems, S\=ugaku {\bf47} (1995), 124--144; English
translation, S\=ugaku Exposition {\bf9} (1997), {\bf??} , A.M.S., Providence

 \bibitem{25} S. Mukai, Curves and symmetric spaces. I, Amer. J. Math.
{\bf117} (1995), 1627--1644

 \bibitem{26} S. Mukai, Curves and K3 surfaces of genus 11, in Moduli of
vector bundles, M. Maruyama (Ed.), Marcel Dekker, New York, 1996,
pp.~189--197

 \bibitem{27} D. Mumford, Projective invariants of projective structures and
applications, in Proc. Internat. Congress of Math., Stockholm, 1962,
pp.~526--530

 \bibitem{28} D. Mumford, Abelian varieties, Oxford Univ. Press, 1970

 \bibitem{29} D. Mumford, J. Fogarty and F. Kirwan, Geometric invariant
theory, Third enlarged edition, Springer, 1994

 \bibitem{30} M. S. Narasimhan and S. Ramanan, Moduli of vector bundles on a
compact Riemann surface, Ann. of Math. {\bf89} (1969), 19--51

 \bibitem{31} M. S. Narasimhan and S. Ramanan, Geometry of Hecke cycles. I, in
C.~P. Ramanujam---a tribute, Springer, 1978, p.~291--345

 \bibitem{32} M. S. Narasimhan and S. Ramanan, $2\theta$-linear systems on
Abelian varieties, in Vector bundles on algebraic varieties (Bombay, 1984),
Tata Inst. of Fundamental Research Studies {\bf11}, Oxford Univ. Press,
1987, pp.~415--427

 \bibitem{33} M. S. Narasimhan and C. S. Seshadri, Stable and unitary vector
bundles on a compact Riemann surface, Ann. of Math. {\bf82} (1965), 540--564

 \bibitem{34} P. Newstead, Characteristic classes of stable bundles, Trans.
Amer. Math. Soc., {\bf169} (1972), 337--345

 \bibitem{35} P. Newstead, Introduction to moduli problems, Tata Inst. lecture
notes, Springer, 1978

 \bibitem{36} P. Newstead, Brill--Noether problem list update, Univ. of
Liverpool, 1992; most recent draft available from the EP VBAC section of
the AGE and EuroProj file server http://euclid.mathematik.uni-kl.de/

 \bibitem{37} D. G. Northcott, Finite free resolutions, Cambridge Tracts in
Math. {\bf71}, Cambridge Univ. Press, 1976

 \bibitem{38} P. Pragacz, Enumerative geometry of degeneracy loci, Ann. Sci.
\'Ecole Norm. Sup. {\bf21} (1988), 413--454

 \bibitem{39} P. Pragacz, Algebro-geometric applications of Schur S- and
Q-polynomials, in Topics in invariant theory, LNM {\bf1178}, Springer
pp.~130--191

 \bibitem{40} S. Ramanan, The moduli space of vector bundles over an algebraic
curve, Math. Ann. {\bf200} (1973), 69--84

 \bibitem{41} C. S. Seshadri, Moduli of vector bundles with parabolic
structures, Bull. Amer. Math. Soc. {\bf83} (1977), 124--126

 \bibitem{42} C. S. Seshadri, Vector bundles on curves, Contemp. Math {\bf153}
(1993), 163--200

 \bibitem{43} M. Thaddeus, Conformal field theory and the cohomology ring
of the moduli space of stable bundles,  J. Diff. Geom. {\bf35} (1992),
131--149

 \bibitem{44} M. Teixidor, Brill--Noether theory for stable vector bundles,
Duke Math. J. {\bf62} (1991), 385--400

 \bibitem{45} J. Wahl, Gaussian maps on algebraic curves, J. Diff. Geom.
{\bf32} (1990), 77--98

 \bibitem{46} J. Wahl, Introduction to Gaussian maps on an algebraic curve,
in Complex projective geometry,  G. Ellingsrud, C. Peskine, G. Sacchiero and
S.A. Stromme (Eds.), LMS lecture notes {\bf179} (1992), pp.~304--323

 \bibitem{47} G. E. Welters, A theorem of Gieseker--Petri type for Prym
varieties, Ann. Sci.  \'Ecole Norm. Sup. {\bf18} (1985), 671--683

 \bibitem{48} D. Zagier, On the cohomology of moduli spaces of rank two
vector bundles over curves, in The moduli space of curves, R. Dijkgraaf and
others (Ed.), Birkh\"auser, 1995, pp.~533--563

\medskip Added in translation:

 \bibitem{Beauville} A. Beauville, Fibr\'es de rang deux sur une courbe,
fibr\'e determinant et fonctions theta. II, Bull. Soc. math. France
{\bf119} (1991) 259--291

 \bibitem{Laszlo} Y. Laszlo, Un th\'eor\`eme de Riemann pour les diviseurs
theta sur les espaces de modules de fibr\'es stable, Duke Math J. {\bf64}
(1991) 333--347

 \bibitem{dual K3s} S. Mukai, Duality of polarized K3 surface, Nagoya
preprint Apr 1997, to appear in Proc.\ of Algebraic geometry Euroconference
(Warwick, 1996), 17~pp.

 \bibitem{PR} K. Paranjape and S. Ramanan, On the canonical ring of a
curve, in Algebraic geometry and commutative algebra (in honour of
M.~Nagata), Vol.~II, Kinokuniya, Tokyo 1988

 \bibitem{OP} W. Oxbury and C. Pauly, Subvarieties of Heisenberg
invariant quartics and $\mathcal{SU}_C(2)$ for a curve of genus four, Duke
file server, alg-geom/9703026, 36~pp.

 \bibitem{OPP} W. Oxbury, C. Pauly and E. Previato, Subvarieties of
$\mathcal{SU}_C(2)$ and $2\theta$ divisors in the Jacobian, Duke file
server, alg-geom/9701010, 41~pp.

 \end{thebibliography}
 \end{document}